\documentclass[11pt]{article}

\usepackage{amssymb}
\usepackage{cite}
\usepackage{epsfig}
\usepackage{rotate}

\newcommand{\beq}{\begin{equation}}
\newcommand{\eeq}{\end{equation}}
\newcommand{\Leq}[1]{\label{#1}\end{equation}}
\newcommand{\bdm}{\begin{displaymath}}
\newcommand{\edm}{\end{displaymath}}

\newtheorem{theorem}{Theorem}

\newtheorem{lemma}[theorem]{Lemma}
\newtheorem{proposition}[theorem]{Proposition}

\def\M{{\cal M}}
\def\O{{\cal O}}
\def\R{\mathbb R}
\def\Q{\mathbb Q}
\def\Z{\mathbb Z}
\def\mod#1{\,({\rm mod\ }#1) }
\def\proof{\medskip\noindent {\bf Proof.}\hskip 10pt}
\def\qed{$\square$}

\title{Geometric representation of interval exchange 
maps over algebraic number fields}
\author{G. Poggiaspalla, J. H. Lowenstein\dag, and F. Vivaldi}
\date
{\it\small
School of Mathematical Sciences, Queen Mary, University of London,
London E1 4NS, UK \\
\dag Dept.~of Physics, New York University, 2 Washington Place, New
York, NY 10003, USA
}

\begin{document}
\renewcommand{\labelenumi}{$(\roman{enumi})$}
\maketitle
\begin{abstract}
This paper is concerned with the restriction of interval 
exchange transformations to algebraic number fields, which 
leads to maps on lattices. 
We characterize renormalizability arithmetically, and study its
relationships with a geometrical quantity that we call the drift vector. 
We exhibit some examples of renormalizable IET
with zero and non-zero drift vector and carry out some
investigations of their properties. In particular we look for evidence
of the finite decomposition property on a family of IETs extending
the example studied in \cite{LPVYocc}.
\end{abstract}
\vspace*{50pt}\centerline{\Large \today}

\section{Introduction}\label{section:Introduction}
An interval exchange transformations (IET) is a one-dimensional map
$E$ of an interval, consisting of a rearrangement of a partition of
that interval. Let $\{\Omega_i\}$ be an $N$-elements partition of
$\Omega=[0,1)$, and let $\pi$ be a permutation of the integers
$1,\ldots,N$. We define $E$ as the map that rearranges the intervals
according to the permutation $\pi$, preserving orientation.
This procedure defines $N$ translation vectors $\tau_i$
\begin{equation}\label{eq:TranslationVectors}
E(x)=x+\tau_i \hskip 30pt x\in\Omega_i.
\end{equation}
We assume that this permutation is {\em irreducible}, that is,
$\pi\{1,\ldots,k\}=\{1,\ldots, k\}$ if and only if $k=N$. 

There is substantial literature on the subject. An early theorem of
Keane \cite{Keane1} states that if the permutation is irreducible
and the set of lengths rationally independent, then two
discontinuity points never lie on the same orbit and $E$ is minimal.
A stronger results was later obtained by Masur and by Veech
\cite{Masur,Veech}: given an irreducible permutation, for almost all
set of lengths, the associated interval exchange map is minimal and
uniquely ergodic.

Boshernitzan \cite{boshernitzan} gave the following ergodicity
criterion: if we let $\epsilon_k$ be the shortest length of an
interval of the partition of $E^k$, then if
$$
\lim_{k\rightarrow\infty}k\epsilon_k>0
$$
$E$ is uniquely ergodic. A result of Katok \cite{Katok} ensures that
interval exchanges cannot be mixing. However, a recent result by
Avila and Forni \cite{avilaforni} shows that, for any permutation
which is not a rotation (meaning that $\pi(i)\not\equiv i+1\mod N$),
the interval exchange map associated to Lebesgue-almost all set of
lengths is weakly-mixing.

Here we consider a uniquely ergodic interval exchange
map  $E$, and we assume that the lengths of the 
sub-intervals $\Omega_i$ generate an algebraic 
extension $K$ of degree $n$ of the rationals.
The {\it module\/} of $E$ is the set of all linear integral
combinations of the sub-interval lengths, namely
\begin{equation}\label{eq:M}
\M=\Bigl\{\sum_{i=1}^N m_i|\Omega_i|\,:\, m_i\in\Z\Bigr\}
\end{equation}
and the {\it rank\/} of $E$ is the rank of its module.
Algebraic IETs constitute a countable, hence exceptional, family.
Little is known about these maps, the only substantial result being
the Boshernitzan and Carroll theorem \cite{BoshernitzanCarroll}:
if an IET is defined over a quadratic field, then, up to rescaling,
the number of induced maps on sub-intervals is finite.
In short, quadratic IETs can be renormalized.
In general, interval exchanges over algebraic number fields
of higher degree do not have this nice property. If they do, 
then their rank must be greater than two.

The code $\omega=\omega(x)$ of a point $x\in\Omega$
is constructed with respect to the natural partition 
$$
\omega_j=i\quad\Leftrightarrow\quad E^{j-1}x\in\Omega_i.
$$
Let $e_1,\ldots,e_N$ be the canonical basis of $\mathbb{Z}^N$. 
We associate to each word $\omega$ its sub-word $\omega^{(k)}$ 
of length $k$, and the sequence
$s(\omega^{(k)})$ in $\mathbb{Z}^N$ defined as follows
\begin{equation}\label{eq:Staircase}
s(\omega^{(k)})=\sum_{j=1}^k e_{\omega_j}
\qquad k=1,2,\ldots.
\end{equation}
We call it the {\it staircase function.}
If we denote by $|\omega^{(k)}|_i$ the number of occurrences of
$i$ in $\omega^{(k)}$, then we have
$$
s(\omega^{(k)})=\sum_{i=1}^N |\omega^{(k)}|_i e_i.
$$
This expression may be rewritten as follows
$$
s(\omega^{(k)})=k \Lambda+D_k.
$$
Here $\Lambda$ is the vector of the lengths of the intervals
\begin{equation}\label{eq:Lambda}
\Lambda=(|\Omega_1|,\ldots,|\Omega_N|)
\end{equation}
and $D_k$ is the vector of the local discrepancies of the sequence
$\omega$ in each cylinder $(i)$ or, equivalently, of the sequence
$(E^k x)$ in each intervals $\Omega_i$
\begin{equation}\label{eq:Discrepancies}
D_k=(D_k(x,\Omega_1),\ldots,D_k(x,\Omega_N))
\end{equation}
and
$$
D_k(x,\Omega_i)=\sum_{j=1}^k {\bf 1}_{\Omega_i}(E^{j=1} x)-k|\Omega_i|
 =|\omega^{(k)}|_i-k|\Omega_i|,
$$
where ${\bf 1}_{\Omega_i}$ is the characteristic function of
the interval $\Omega_i$.
By Birkhoff's ergodic theorem and using minimality we have
$||D_k||=o(k)$; the staircase of each point thus has the 
principal direction $\Lambda$ with sub-linear perturbations.

The theory of the deviation from the ergodic averages in IET is closely
related to the theory of geodesic flows on flat-surfaces, for which
strong results have been given by Zorich \cite{zorichdev,zorichgaussm}.
Specifically, we suppose that the IET can be described as a section
of a geodesic flow on a flat surface of genus $g$.
Then, for almost all such IETs, there is a nested sequence of subspaces
$$
\mathcal{H}_1\subset\mathcal{H}_2\subset\cdots\subset\mathcal{H}_g\subset\R^N
$$
with $\mathcal{H}_1$ being the one-dimensional subspace spanned by
$\Lambda$, and each $\mathcal{H}_i$ being $i$-dimensional, such that,
for any $f\in\mathcal{H}_i^\perp$ with $||f||=1$
$$
\limsup_{k\to\infty}\frac{\log\langle s(\omega^{(k)}),f\rangle}{\log k}=\nu_i.
$$
(The symbol $\langle\,,\rangle$ denotes the standard inner product.)
This result says that the deviations of the staircase function from the
ergodic means follow some power laws, the exponents of which are
given by the Lyapunov exponents $\nu_i$, depending only on the
so-called Rauzy class the IET belongs to. 

The aim of this paper to introduce a lattice representation of 
algebraic IETs, following a method developed for higher dimensional 
piecewise isometries 
\cite{KouptsovLowensteinVivaldi,LowensteinKouptsovVivaldi},
and analogous to the representation used in \cite{LPVYocc}.
Since the translations $\tau_i$ are expressed as integral linear 
combinations of the lengths, for any point $\xi$ in $K\cap [0,1)$
the set $\xi+\M$ is invariant under $E$ and contains the orbit of $\xi$. 
We then identify $\xi+\M$ with the lattice ${\bf L}_\xi\cong\mathbb{Z}^m$,
where $m\leq n$ is the rank of $\M$.
This construction may be seen as a projection of the lattice
$\Z^N$ in which the staircase is built onto the lattice ${\bf L}$
where the dynamics takes place. In this way we obtain a discrete
geometric representation of the IET, restricted to a dense set of 
algebraic points.

The content of this paper is the following:

In section \ref{latticendrift} we introduce a module-theoretic
description of an algebraic IET and define a fundamental
dynamical invariant, the {\em drift vector} $\mathcal{S}$.
We then show that the drift vector is a representation of 
the classical Sah-Arnoux-Fathi invariant (proposition 
\ref{proposition:SAF}). Most of the constructs developed
in this section refer to the case of maximal rank;
this restriction is justified by the fact that renormalizable 
IETs have this property (see theorem \ref{theorem:Main} below).

In section \ref{renorm} we deal with renormalizability of 
algebraic IETs and we describe the role of the drift vector in the
renormalization process. We then prove the main result of 
this section, which appears as theorem \ref{theorem:rhounit} and 
proposition \ref{proposition:DriftConditions}.
(For background on substitution dynamical systems see \cite{PytheasFogg}.) 
\begin{theorem}\label{theorem:Main}
Let $E$ be a uniquely ergodic interval exchange transformation on 
$N$ intervals, which is renormalizable with scaling ratio $\rho$.
Then $\rho$ is a unit in the ring $\Z[\beta]$,
where $\beta=\rho^{-1}$ is the largest eigenvalue of the incidence 
matrix of the primitive substitution associated with $E$. 
The map $E$ is defined over $\mathbb{Q}(\rho)$, its
rank $n$ is maximal and equal to the degree of $\rho$, and 
\begin{enumerate}
\item {if the drift vector is non-zero, then $n$ is even and $n\leq N$;}
\item {if the drift vector is zero, then $2n\leq N.$}
\end{enumerate}
\end{theorem}

Subsequently, we introduce the recursive tiling property and
formalism (Vershik coding), and give estimates for the
asymptotic behaviour of points with eventually periodic code 
in the scaling dynamics (propositions \ref{proposition:Bound}
and \ref{proposition:AsymptoticPeriodic}).
We also give a condition for these points to be of positive 
density on the lattice (proposition \ref{proposition:PositiveDensity}).

Section \ref{examples} is devoted to the analysis of specific examples.
Using Rauzy induction, we construct a renormalizable IET of 
degree four over four intervals, with non-zero drift vector. 
We then construct an infinite family of cubic renormalizable IET 
with zero drift-vector, which belong to the so-called Arnoux-Rauzy 
family \cite{ArnouxRauzy}. 
We also perform numerical investigations which give strong 
evidences that some of these examples share the {\em finite 
decomposition property} encountered in \cite{LPVYocc}.
This property, which we establish in one case (proposition 
\ref{proposition:E2}) is described as follows. 
The field $K$ may be represented as the disjoint union of 
invariant layers $\xi+\M$, for a suitable set of representative 
points $\xi\in K$. 
An algebraic IET has the finite decomposition property if
the restriction of each layer to the unit interval is the union 
of finitely many orbits.

These results show connections between the finite decomposition 
property of a renormalizable IET and the Pisot property of the 
(reciprocal of the) scaling constant. 
This problem deserves further investigation.

\bigskip

\noindent {\sl Acknowledgements.\/} This research was supported by
EPSRC grant No GR/S62802/01. One of the authors (GP) wants to
thank S. Ferenczi, X. Bressaud and P. Arnoux for some enlightening
discussions.

\section{Lattice dynamics and drift vector}\label{latticendrift}

In this section we consider algebraic IETs with maximal rank,
for which an interesting arithmetical theory arises. 
In the next section we shall see that renormalizable IETs 
have this property.

\subsection{Lattices}\label{section:Lattices}
We choose a basis for $\mathcal{M}$ as a $\mathbb{Z}$-module,
and consider the translation vectors $\tau_i\in\mathcal{M}$
defined by $E$ ---cf. equation (\ref{eq:TranslationVectors}).
Let $\phi$ be the map from
$\mathcal{M}$ in ${\bf L}$ sending the basis elements of
$\mathcal{M}$ onto the canonical basis of ${\bf L}$.
We have two lattices, one is the symbolic
lattice $\mathbb{Z}^N$ in which the staircase was defined; the other
is ${\bf L}$, it has dimension $n$ and represents the module in
which the dynamics takes place. We call $\pi$ the projection from
the former in the latter; this is a linear map such that
\begin{equation} \label{eq:vi}
\pi(e_i)=v_i \hskip 40pt 
v_i:=\phi(\tau_i)
\end{equation}
which is represented by an $n\times N$ matrix.
We denote by $\psi=\phi\circ E \circ \phi^{-1}$ the map
conjugate to $E$ on the lattice {\bf L}.

As we have a map of the unit interval, the dynamics of $\psi$ does
not take place on all of ${\bf L}$, but only on the invariant slab
of points $z\in{\bf L}$ such that $\phi^{-1}(z)\in [0,1)$.
We are then led to consider the set 
$\mathcal{M}\cap [0,1)$
for which we require a suitable representation.
Let
\begin{equation}\label{eq:O}
\O=\{\zeta\in K: \zeta\M\subset\M\}
\end{equation}
where $K$ is the smallest field containing $\mathcal{M}$.
Then $\O$ is an order\footnote{A finitely-generated subring of $K$, which
contains 1 and has maximal rank.} in $K$. 
Let $\nu_1,\ldots,\nu_{n}$ be a $\Z$-basis for $\M$, and let $d$ be 
the smallest positive integer such that $d\nu_k\in\O$, for $k=1,\ldots,n$.
Then $J=d\M$ is an ideal of $\O$ \cite[section 4.7]{Cohen};
the normalized generators $d\nu_k,\,k=1,\ldots,n$, are algebraic
integers which form a $\Z$-basis for $J$.
By construction, $\M$ is also a module over $\O$; 
furthermore $\O$ is the largest subring of $K$ with this property.
Now $\sum_{i=1}^N|\Omega_i|=1\in\M$, and so $\Z\subset \M$.
We can form the quotient
$$
\M':=\M/\Z
$$
which is a $\Z$-module, which we shall represent as $\M\cap[0,1)$. 
Let $j$ be the smallest positive integer in $J$: then the rational numbers
in $\M$ are the multiples of $j/d$, and hence
\begin{equation}\label{eq:Torsion}
(\Q\cap\M)/\Z=\frac{j}{d}\,\Z/\Z\cong\Z/b\Z \hskip 30pt
b=\frac{d}{\mbox{gcd}(d,j)},
\end{equation}
that is, there are $b$ congruence classes of rational numbers in
$\M'$, while the non-rational part is freely generated. Thus we have
a $\Z$-module isomorphism $\M' \to {\bf L'}$ where
\begin{equation}\label{eq:Isomorphism}
{\bf L'}=(\Z/b\Z)\oplus\Z^{n-1}
\end{equation}
which characterizes the additive structure associated with the
map $E$. Accordingly, we choose for $\zeta\in\M'$ the following
representation
\begin{equation}\label{eq:Zeta}
\zeta=\frac{1}{d}\sum_{k=1}^{n}m_k\,\nu'_k \hskip 30pt
m_k\in\Z\qquad \nu'_k\in J,\quad \nu'_1=j,
\end{equation}
where the $\nu'_k$ are obtained from the normalized generators
$d\nu_k$ via a unimodular transformation. The isomorphism
(\ref{eq:Isomorphism}) then becomes
\begin{equation}\label{eq:IsomorphismII}
\phi': \,\zeta\mapsto(m_0\mod{b},m_1,\ldots,m_{n-1}).
\end{equation}


To represent the dynamics in the field $K$ ---the extension of the
rational generated by the intervals' lengths--- we use the
$\Z$-basis $(\nu_k)$ for $\mathcal{M}$, and any point
$x$ in $K$ can be written as
$$
x=\sum_{k=1}^n r_k\nu_k\quad\textrm{for}\quad r_k\in\Q.
$$
By definition, $E$ preserves $\mathcal{M}$ and for any
$x\in K$ its orbit is confined in $x+\mathcal{M}$. We have
a ``layered'' representation of $K$ and we regard it as the union of
disjoint invariant sets of the form
\begin{equation}\label{layers}
K=\bigcup_{\xi\in\Xi}(\xi+\mathcal{M})
\end{equation}
for a chosen set of representative $\Xi$.
Following \cite{LPVYocc}, we choose
\begin{equation}\label{eq:Xi}
\Xi=\Bigl\{\sum_{k=1}^n r_k\nu_k,\quad r_k\in\Q\cap [0,1)\, \Bigr\}.
\end{equation}
The dynamics of $E$ acts on each layer separately, which in 
each case may be represented as the lattice ${\bf L}$, with a 
$\xi$-dependent dynamics on it.

In general, the closer two points get in $\mathcal{M}$, the
further away their images in the lattice ${\bf L}$ will be.
Indeed, let $x\in\mathcal{M}$ correspond to the lattice point
$z=(z_1,\ldots,z_n)$
\begin{equation}\label{eq:xz}
x=\sum_{k=1}^n z_k\nu_k\in\mathcal{M}.
\end{equation}
If the coefficients $z_k$ are bounded, then there are only finitely
many possibilities for $x$, and hence $x$ is bounded away from zero.
Thus if $|x|\to 0$ then $\Vert z\Vert\to\infty$, and one is 
interested in the relation between these limits.

We recall that height of a polynomial is the maximum absolute 
value of the coefficients, and the height of an algebraic number 
is the height of its minimal polynomial.
In the particular case when $\mathcal{M}$ takes the form
$\M=\Z[\lambda]$ for some algebraic number $\lambda$ of degree $n$,\
\footnote{$\Z[\lambda]$ is the $\Z$-module
generated by $1,\lambda,\lambda^2,\ldots,\lambda^{n-1}$}
Liouville's inequality gives the following bound.
\begin{proposition}\label{proposition:Liouville}
If $\mathcal{M}=\Z[\lambda]$ for some algebraic number 
$\lambda$ of degree $n$, then for every $x\in\M$ we have
\begin{equation}\label{eq:Liouville}
|x|\,\Vert z\Vert^c\geq e^{c(n-1)}
\hskip 40pt
z=\phi(x)\quad c=\deg\lambda+\log H(\lambda)
\end{equation}
where $H(\lambda)$ is the height of $\lambda$.
\end{proposition}
{\bf Proof:} We have (cf. \cite{Roy}) that for an integer polynomial
$f$ with height $\leq H$ and degree $d$, for all algebraic number
$\lambda$ we have either
$$
|f(\lambda)|\geq e^{-c(d+\log H)}
$$
or $\lambda$ is a root of $f$. 
The constant is $c=\deg\lambda+H(\lambda)$.
In our context, the polynomial $f$ has coefficients $z_i$, and
$\lambda$ is the generator of $\mathcal{M}$, that is, the
basis vectors $\nu_k=\lambda^{k-1}$.
Then $H=||z||$ for the max-norm and
$$
|f(\lambda)|=|x|\geq e^{-c(n-1)} ||z||^{-c},
$$
as desired.
\qed

\subsection{Drift vector}

Let $v_i$ be the image on the lattice of the translation vector
$\tau_i$ ---see equation (\ref{eq:vi}). The {\it drift vector\/}
$\mathcal{S}$, defined as
\begin{equation}
\mathcal{S}=\sum_{i=1}^N |\Omega_i|v_i=\pi(\Lambda)\label{Sdef}
\end{equation}
gives the dominant direction of the dynamics over the lattice.
The drift vector is related to a known invariant of IET theory,
the so-called {Sah-Arnoux-Fathi invariant,} as we shall see in
section \ref{section:saf}.

Let the point $x\in\M$ have itinerary $\omega$, and let
$z=\phi(x)$ be the corresponding lattice point.
Then the orbit of $z$ may be expressed in terms of the
drift vector as follows
\begin{equation}\label{eq:Alphak}
z+\pi s(\omega^{(k)})=z+k \mathcal{S}+\pi D_k\nonumber
                     =z+k(\mathcal{S}+\alpha(k))
\end{equation}
where the staircase function $s$ and the discrepancies $D_k$ were
defined in equations (\ref{eq:Staircase}) and (\ref{eq:Discrepancies}),
respectively, and $\alpha(k)=\pi D_k/k$.
Because $|D_k|=o(k)$, the orbits of any lattice map with nonzero
drift vector will be dominated by a linear drift in the direction
of $\mathcal{S}$.

\subsection{Density}

The dynamics in ${\bf L}$ takes place on a subset of zero density.
We can nonetheless have a meaningful measure of the ``size'' of the
orbits by defining a modified density function on the reduced
lattice ${\bf L'}$, which results from restriction to the unit interval. 
With reference to (\ref{eq:Torsion}) and
(\ref{eq:Isomorphism}), we represent a subset $A$ of $\M'$ as
$$
A=A_0+\cdots+A_{b-1}\hskip 30pt A_i\cap \sum_{k\ne i}A_k=\{0\}
\hskip 30pt A_i\cap \Q\equiv \frac{i}{b}\mod{\Z}.
$$
Letting ${\bf A}_j=\phi(A_j)$, we define
$$
D({\bf A}_j)=\lim_{k\rightarrow\infty}\frac{1}{b (2k)^{n-1}}
\sharp\{(m_0\mod{b},m_1,\ldots,m_{n-1})\in {\bf A}_j,\,|m_i|\leq
k,i\neq 0\}.
$$
Note that $m_0\mod{ b}$ is fixed, by construction. If the 
limit above exists for all $j$, we define the density of ${\bf A}$ as
$$
\mathcal{D}({\bf A})=\sum_{j=0}^{b-1} \mathcal{D}({\bf A}_j)
$$
which is a (non countably additive) probability measure on ${\bf
L'}$.

The density $\mathcal{D}$ is invariant under an unimodular
transformation $U$. The density $\mathcal{D}$ measures the
proportion of the points lying inside an hyper-cube 
$\mathcal{C}_k$ centred at the origin, with side $2k$,
as $k$ becomes large.
Considering $\mathcal{D}(U{\bf A})$ is equivalent to consider
$\mathcal{D}'({\bf A})$ where $\mathcal{D}'$ measures the limiting
proportion of the points lying inside the rhombus $U^{-1}\mathcal{C}_k$,
also centred on the origin.
As for all $k$ there exists $k'$ and $k''$
such that $U^{-1}\mathcal{C}_{k'}\subset\mathcal{C}_k\subset
U^{-1}\mathcal{C}_{k''}$ we deduce that $\mathcal{D}({\bf
A})=\mathcal{D}'({\bf A})=\mathcal{D}(U{\bf A})$.

\begin{proposition}\label{proposition:Density}
Let $I\subset [0,1)$ be a finite union of intervals, and let 
${\bf A}=\phi(\M\cap I)$. Then $\mathcal{D}({\bf A})$ exists 
and is equal to the Lebesgue measure of $I$.
\end{proposition}
\proof We begin with the case in which $I$ is a single interval.
Let ${\bf m}=(m_0,\ldots,m_{n-1})$, let 
$x_{\bf m}=\phi^{-1}({\bf m})$, and let
$$
C_k=\frac{1}{b(2k)^{n-1}}\sharp\{{\bf m}\in{\bf A}
         \,:\,|m_i|\leq k\}\qquad k=1,2,\ldots.
$$
Now let $m_0,m_2,\ldots,m_{n-1}$ be fixed integers,
and define
$$
C^\prime_{m_0,k,m_2,\ldots,m_{n-1}}=\frac{1}{2k}\sharp
\{(m_0,m,m_2,\ldots,m_{n-1})\,:\, |m|\leq
k,~x_{m_0,m,m_2,\ldots,m_{n-1}}\in I\}
$$
then by Weyl's criterion
$$
C_{m_0,k,m_2,\ldots,m_{n-1}}-|I|
\,\rightarrow\,
0\quad\textrm{when}\quad k\,\rightarrow\,\infty.
$$
We then have
$$
C_k=\frac{1}{b(2k)^{n-2}}\sum_{m_0=0}^{b-1}\sum_{m_2,\ldots,m_{n-1}=-k}^k
C^\prime_{m_0,k,m_2,\ldots,m_{n-1}}
$$
and
$$
C_k=\frac{1}{b(2k)^{n-2}}\sum_{m_0=0}^{b-1}\sum_{m_2,\ldots,m_{n-1}=-k}^k
\left(\epsilon_{m_0,k,m_2,\ldots,m_{n-1}}+|I|\right).
$$
where $\epsilon_{m_0,k,m_2,\ldots,m_{n-1}}=
C_{m_0,k,m_2,\ldots,m_{n-1}}-|I|$.
Finally
$$
C_k=|I|+\frac{1}{b(2k)^{n-2}}\sum_{m_0=0}^{b-1}\sum_{m_2,\ldots,m_{n-1}=-k}^k
\epsilon_{m_0,k,m_2,\ldots,m_{n-1}}
$$
the last term is a $(n-2)$-dimensional Ces\`aro sum and thus converges
to zero when $k\rightarrow\infty$, leading to
$\mathcal{D}({\bf A})=|I|$. 
The argument above is easily extended to
the case of a finite union of intervals.
\qed

Next we have
\begin{proposition}\label{proposition:NonZeroDrift}
If the drift vector is non-zero and if $n>2$, then the density 
of any orbit on the lattice ${\bf L'}$ is zero.
\end{proposition}

\proof This is an easy consequence of formula (\ref{eq:Alphak}).
We have
$$
||\psi^k(z)-z||=||\mathcal{S} k+k\alpha(k)||\geq
k\Big|||\mathcal{S}||-||\alpha(k)||\Big|\geq k\min_k
\Big|||\mathcal{S}||-||\alpha(k)||\Big|.
$$
This latter minimum is attained because $1>||\alpha(k)||\rightarrow
0$. Thus the number of points inside the box centred at the origin
and with side $2k$ increases at most linearly with $k$. So, if the 
degree of the IET is greater than 2, then the lattice {\bf L}' has 
rank greater than one, leading to zero density.
\qed

In particular, an interval exchange transformation can have the
{\it finite decomposition property\/} ---each of its layers is the 
union of a finite number of orbits--- only if it has zero drift. 
The actual value of the drift vector depends on the choice 
of the basis of $\M$; changing the basis is equivalent to apply 
a unimodular transformation on $\mathcal{S}$ which would indeed 
change its direction and magnitude. Nevertheless, we will see that
the relevant property of $\mathcal{S}$ is mainly its nullity.

\subsection{The Sah-Arnoux-Fathi invariant}\label{section:saf}
In this section we regard $\mathbb{R}$ as a vector space over 
$\mathbb{Q}$, and consider the tensor product 
$\mathbb{R}\otimes\mathbb{R}$ of $\mathbb{Q}$-vector spaces.
To a general $N$-interval exchange $E$ with set of length
$|\Omega_i|$ and translation vectors $\tau_i$, Arnoux
\cite{thesearnoux} associates an element $S$ of $\R\otimes\R$. The
quantity $S$ can be computed by the formula
$$
S=\sum_{i=1}^{N} |\Omega_i|\otimes\tau_i.
$$
The Sah-Arnoux-Fathi (SAF) invariant is a non-trivial conjugacy invariant; 
it is not a complete invariant but is in some sense the best possible one. 
The set of all the IETs constitute a group under the composition law. 
The map associating to each IET its invariant $S$ is a homomorphism 
from the group of IETs into $\mathbb{R}\otimes\mathbb{R}$, with the 
smallest possible kernel. 
For further details and discussions about the algebraic
significance of $S$ see \cite{thesearnoux,veech3}. 
From these references we also know that $S$ is an induction invariant 
for minimal transformations, i.e., $S$ is the same computed on any
induced maps of a given minimal interval exchange transformation. 
If $E$ is periodic, then its invariant is zero. The converse is
true for an interval exchange of three sub-intervals, but not
for more subintervals, as shown in \cite{thesearnoux}.
The Arnoux-Yoccoz maps \cite{Arnoux,ArnouxYoccoz} provide an example
of this phenomenon: they have 7 sub-intervals and zero invariant,
but they are not periodic, in fact they are minimal and uniquely ergodic.

In general, $\R\otimes\R$ is infinite dimensional, but in our case
the lengths and the translations belong to $K$, which is
a vector space over $\Q$ of dimension $n$. In fact we can even restrict
ourselves to the $\Z$-module $\M$. The quantity $S$ should then be
in $\M\otimes_\Z\M$ (the subscript $\mathbb{Z}$ emphasizes
that the tensorial product is taken between $\mathbb{Z}$-modules).
As above, we fix a $\Z$-basis $\nu_1,\ldots,\nu_{n}$ for $\M$. 
We define a linear map $V: \M\otimes_\Z\M \to \M^{n}$ by
specifying its action on this basis as follows
$$
V: \nu_i\otimes\nu_j\longmapsto
\underbrace{(0,\ldots,\nu_i,\ldots,0)}_{
\textrm{ $\nu_i$ has position $j$}}.
$$
The map $V$ is bijective and is an isomorphism of
$\mathbb{Z}$-modules. We have
$$
V(S)=\sum_{i=1}^{N}V(|\Omega_i|\otimes\tau_i).
$$
Next we express the translation vectors $\tau_i$ and 
lengths $|\Omega_i|$ with respect to the module basis
$$
\tau_i=\sum_{k=1}^{n}\tau_i^{(k)}\nu_k
\hskip 40pt
|\Omega_i|=\sum_{j=1}^{n}|\Omega_i|^{(j)}\nu_j
$$
so that
$$
|\Omega_i|\otimes\tau_i=\sum_{j,k=1}^{n}|\Omega_i|^{(j)}
\tau_i^{(k)}\nu_j\otimes\nu_k.
$$
Finally, recalling that $v_i=\phi(\tau_i)$ are the images of
the translation vector on ${\bf L}$ (equation (\ref{eq:vi})), we have
\begin{eqnarray*}
V(S)&=&\sum_{i=1}^{N}
 (\tau_i^{(1)}\sum_{j=1}^{n}|\Omega_i|^{(j)}\nu_j,\ldots,\tau_i^{(n)}
 \sum_{j=1}^{n}|\Omega_i|^{(j)}\nu_j)\\
&=&\sum_{i=1}^{N}(\tau_i^{(1)},\ldots,\tau_i^{(n)})|\Omega_i|
 =\sum_{i=1}^{N}v_i|\Omega_i|=\mathcal{S}
\end{eqnarray*}
which shows that the drift vector $\mathcal{S}$ is a representation 
of the SAF invariant in our specific context.

The condition $\mathcal{S}=0$ is quite strong; in fact
it implies $n$ rational dependencies between the lengths of the
intervals. Note that if we have only two intervals, the two lengths
are of the form $\Omega_1$ and $1-\Omega_1$ and the condition
$\mathcal{S}=0$ would imply $|\Omega_1| k+(1-|\Omega_1|)l=0$ for two
integers $k$ and $l$, which is impossible since $|\Omega_1|$ is
irrational. In fact a two-interval exchange is a pure rotation of
the circle, which gives a pure translation on the lattice. More
generally, for $N$ intervals, we have $N$ lengths $|\Omega_i|\in\M$,
which we assumed to have rank $n$. For $\mathcal{S}$ to be zero, we
must have
$$
\sum_{i=1}^{N}x_i|\Omega_i|=\sum_{i=1}^{N-1}(x_i-x_{N})|\Omega_i|+x_{N}=0
$$
for some integers $x_i$. Then a necessary condition for
$\mathcal{S}=0$ is $n<N-1$. 

The following result summarizes the considerations above.

\begin{proposition}\label{proposition:SAF}
The SAF invariant of an algebraic IET is zero if and only if its
drift vector is zero. If the degree $n$ and the number of intervals
$N$ are such that $n\geq N-1$, then the drift vector, and hence 
the SAF invariant, are non-zero.
\end{proposition}

By the discussion above, the Arnoux-Yoccoz
\cite{Arnoux,ArnouxYoccoz} and Arnoux-Rauzy \cite{ArnouxRauzy}
families of IET have zero SAF
invariant and thus zero drift vector.


\section{Renormalizable cases}\label{renorm}

\subsection{Algebraic consequences of renormalizability}
\label{section:algconsrenorm}

In the following, we assume $E$ to be {\em renormalizable}. 
An IET of $N$ intervals is renormalizable if it has the
same dynamics at each scale level, meaning that there exists an 
interval $I\subset[0,1)$ with the property that the induced map 
on $I$ is conjugated to the IET itself via a simple rescaling $\rho$. 
Thus the induced map on $I$ is an IET of $N$ intervals whose lengths 
are $\rho|\Omega_i|$ with a common scaling factor $\rho$. 
In this section, we prove the following theorem.

\begin{theorem} \label{theorem:rhounit}
Let $E$ be uniquely ergodic interval exchange transformation 
on $N$ intervals, which is renormalizable with scaling ratio $\rho$.
Then $\rho$ is a unit of degree $\leq N$ in the ring $\Z[\beta]$,
where $\beta=\rho^{-1}$ is the largest eigenvalue of the incidence 
matrix of the primitive substitution associated with $E$. 
Furthermore, $E$ is defined over $\mathbb{Q}(\rho)$ and
it has maximal rank.
\end{theorem}

\proof
By construction, the scaled lengths can be expressed as integral 
linear combinations of the lengths $|\Omega_i|$. If we represent 
these linear combinations as an $N$ by $N$ integral matrix $B$, 
then we see that, by definition, the lengths vector $\Lambda$ 
is an eigenvector of $B$ with eigenvalue $\rho$. This tells us 
that $\rho$ is an algebraic integer, and if $n$ is its degree,
then $n\leq N$.
Since $B\cdot\Lambda=\rho\Lambda$ and $\Lambda$ is a probability
vector, then
$$
\sum_{i=1}^{N}(B\cdot\Lambda)_i=\rho\sum_{i=1}^N|\Omega_i|=\rho.
$$
We see that $\rho\mathcal{M}\subset\mathcal{M}$ and
$\rho\in\mathcal{M}$. Now, $\mathcal{M}$ contains 1, and
repeated scaling show that it also contains all positive 
powers of $\rho$. Hence $\M$ contains the ring $\Z[\rho]$. 
It follows that the rank of $\M$ is at least as large as 
the degree $n$ of $\rho$.

A renormalizable and uniquely ergodic interval exchange is 
conjugated to a primitive substitution $\sigma$.
Furthermore, the atoms of the induced map on $I$ have their
iterates cover the whole interval $\Omega$, lest there would be
either some periodic cells or some different ergodic component.
We conclude that
$$
\sum_{i=1}^N\rho^k |\sigma^k(i)||\Omega_i|=1\qquad k\geq 0.
$$
Because the substitution is primitive, for all $i$ the quantity
$|\sigma^k(i)|$ grows like $\beta^k$, the powers of the Perron
eigenvalue of the incidence matrix $M_\sigma$ of $\sigma$, and we have
$$
\rho=\beta^{-1}.
$$
Because both $\rho$ and $\beta$ are algebraic integers, they
are units.

Next we show that $\Lambda$ is the positive eigenvector of $M_\sigma$
associated with $\beta$. Considering the fixed point $\omega$ of
$\sigma$, let $j$ be the unique integer with the property 
that $\sigma^k(j)$ is an approximation to $\omega$.
We find
$$ M_\sigma\left(\begin{array}{c}
|\Omega_1|\\
\vdots\\
|\Omega_N|\end{array}\right)=\lim_{k\to\infty}\frac{1}{|\sigma^k(j)|}M_\sigma\left(\begin{array}{c}
|\sigma^k(j)|_1\\
\vdots\\
|\sigma^k(j)|_N\end{array}\right)=\lim_{k\to\infty}\frac{1}{|\sigma^k(j)|}\left(\begin{array}{c}
|\sigma^{k+1}(j)|_1\\
\vdots\\
|\sigma^{k+1}(j)|_N\end{array}\right).
$$
For all $i$, we have
$$
\frac{1}{|\sigma^k(j)|}|\sigma^{k+1}(j)|_i=
\frac{|\sigma^{k+1}(j)|}{|\sigma^k(j)|}\times
\frac{|\sigma^{k+1}(j)|_i}{|\sigma^{k+1}(j)|}\,\to\,\beta|\Omega_i|.
$$

Since $M_\sigma\Lambda=\beta\Lambda$, for the same reason as
above we have $\beta\in\mathcal{M}$, and as
$\beta\mathcal{M}\subset\mathcal{M}$ and
$\beta^{-1}\mathcal{M}=\rho\mathcal{M}\subset\mathcal{M}$,
we have $\beta\mathcal{M}=\rho\mathcal{M}=\mathcal{M}$.
We see that the vector of the lengths is the solution
of the system of equations
$$
M_\sigma z=\beta z\quad\textrm{and}\quad\sum_{i=1}^N z_i=1
\hskip 40pt z=(z_1,\ldots,z_N).
$$
This solution is unique. Indeed, by the Perron-Frobenius
theorem, $\dim\ker(M_\sigma-\beta {\bf 1})=1$, and the
solution belongs to $\mathbb{Q}(\beta)=\mathbb{Q}(\rho)$. 
Then $\mathcal{M}\subset\mathbb{Q}(\rho)$, and therefore 
the rank of $\M$ is at most $n$. But then the rank is 
exactly $n$, namely, it's maximal.\quad 
$\Box$

The proof above shows that $\mathcal{M}$ is a 
$\Z[\rho]$-module, and hence $\Z[\rho]$ is a subring of the 
order $\O$ ---see equation (\ref{eq:O}) and following remarks.

Let $R$ be the matrix representing the action of multiplication
by $\rho$ on the lattice ${\bf L}$. Since $\rho$ is a unit,
the matrix $R$ is unimodular; moreover, from linear algebra
and the fact that $\rho$ has degree $n$, we have that the
characteristic polynomial of $R$ is equal to the minimal
polynomial of $\rho$.
In particular, $\rho$ is an eigenvalue of $R$.

The codes of the orbits from a scale to the next can be found by
mean of the substitution $\sigma$. If a point $x$ belongs to the
atom $i$, we have
\begin{equation}
E^{|\sigma(i)|}(\rho x)=\rho E(x),\label{scalingconj}
\end{equation}
or on the lattice
\begin{equation}\label{renormpsi}
\psi^{|\sigma(i)|}(R z)=R \psi(z).
\end{equation}
Letting the projection $\pi$ be as in (\ref{eq:vi}), we have
\begin{eqnarray*}
\psi^{|\sigma(i)|}(R z)&=&Rz+\sum_{j=1}^{|\sigma(i)|}v_{\sigma(i)_j}
\,=\,Rz+\sum_{j=1}^{N}|\sigma(i)|_jv_j
\,=\,Rz+\pi\sum_{j=1}^{N}|\sigma(i)|_je_j\\
&=&Rz+\pi M_\sigma\sum_{j=1}^{N}|i|_j e_j=Rz+\pi M_\sigma e_i.
\end{eqnarray*}
By equation (\ref{renormpsi}), we then have
$$
R\psi(z)=Rz+R v_i=Rz+\pi M_\sigma e_i
$$
and since $\pi e_i=v_i$
$$
R\pi e_i=\pi M_\sigma e_i\qquad i=1,\ldots, N.
$$
Finally
\begin{lemma}\label{rpipim}
The following diagram commutes
$$
\matrix{\Z^N & \smash{\mathop{\longrightarrow}\limits^{M_\sigma}} &\Z^N\cr
        \big\downarrow\rlap{$\pi$} &&\big\downarrow\rlap{$\pi$} \cr
        {\bf L}&\smash{\mathop{\longrightarrow}\limits^{R}} &{\bf L}\cr}
$$
\end{lemma}
As a corollary, we see that if $X$ is an eigenvector of $M_\sigma$
with eigenvalue $\zeta$ and $X$ is not in the kernel of $\pi$, 
then $\pi X$ is an eigenvector of $R$ with the same eigenvalue
$$
\pi M_\sigma X=\zeta \pi X=R\pi X.
$$
In addition, the matrix $M_\sigma$ preserves $\ker \pi$. Indeed
if $X\in\ker\pi$, then again $R\pi X=0=\pi M_\sigma X$ and $M_\sigma
X\in\ker\pi$. In particular, the characteristic polynomial of $M_\sigma$
has a factor of degree $N-n$.

We already know that $\Lambda$ is an eigenvector of $M_\sigma$
with the Perron eigenvalue $\beta$.
We thus have two possibilities: $i)$\/ $\pi\Lambda\neq 0$
(the drift vector is non-zero) and $\pi\Lambda$ is an eigenvector
of $R$ with eigenvalue $\beta$;
\/ $ii)$\/ $\Lambda\in\ker\pi$ and $\beta$ is not
part of the spectrum of $R$. This is because the spectrum of
$R$ is contained in the spectrum of $M_\sigma$ and the Perron
eigenvalue is simple.

\subsection{Connections with drift vector}
\label{subsection:ConnectionsWithDriftVector}

When the drift vector $\mathcal{S}$ is non-zero, the algebraic
number $\beta$ is in the spectrum of $R$.
This means that $\rho$ and $\beta=\rho^{-1}$ are algebraic
conjugates, namely roots of the same irreducible polynomial
$p(x)$ of degree $n$.
In fact, all eigenvalues of $R$ are algebraic conjugates
and all are units.

Let now $\tilde p(x)=x^{\deg p}p(1/x)$ be the reciprocal polynomial
of $p$ ---the polynomial whose roots are the reciprocal of the roots
of $p$. Because $p$ is irreducible, so is $\tilde p$, and since
$\beta$ is a root of both $p$ and $\tilde p$, it follows that 
$p=\tilde p$. So the roots of $p$ regroup into reciprocal pairs,
and since these roots cannot include $\pm 1$, the degree $n$ of 
$p$ must be even. 
Any self-reciprocal polynomial $p$ of even degree 
takes the form (see for instance \cite{ahmadi})
$$
p(x)=x^{n/2} h(x+x^{-1})
$$
for some an integral polynomial $h$ of degree $n/2$. The polynomial
$h$ is irreducible because $p$ is. So if $\zeta$ is root of $p$,
then $\xi=\zeta+\zeta^{-1}$ is root of $h$, and since
$\zeta^2-\xi\zeta+1=0$, we see that $\zeta$ lies in a quadratic
extension of the algebraic number field generated by $\xi$.

\def\oldstuff{
{\bf Remark:} Indeed, suppose that the characteristic polynomial of
a matrix $A$ splits into $\mathbb{C}$ as
$$
P_A(t)=\pm (t-\mu_1)^{r_1}\ldots(t-\mu_k)^{r_k}
$$
where $\mu_i$ are the distinct eigenvalues of $A$ and $r_i$ their
multiplicities. Then it is known (cf. for instance \cite{dummit})
that the minimal nullifying polynomial $m_A(t)$ of $A$ has the form
$$
m_A(t)=\pm (t-\mu_1)^{s_1}\ldots(t-\mu_k)^{s_k}
$$
for some integers $s_i\leq r_i$ and $A$ is diagonalizable if and
only if each $s_i=1$.\\
}

Now $p$ is the characteristic polynomial of both $R$ and
$R^{-1}$, which have the same spectrum.
Because $p$ is a minimal polynomial, its roots all have multiplicity 
one and therefore $R$ and $R^{-1}$ are diagonalizable
(see, for instance, \cite[\S 3.10]{jacobson}).
It then follows that $R$ and $R^{-1}$ are conjugate, indeed 
integrally conjugate. That is, there exists an integral matrix 
$G$ such that
$$
R^{-1}=GRG^{-1}
$$
showing that $R$ has time-reversal symmetry \cite{LambRoberts}.
Then $G$ is an involution, exchanging two subspaces of dimension $n/2$.

\def\junkII{
This means that if $v$ is an eigenvector of $R$ with eigenvalue
$\zeta$ then
$$
GRv=G\zeta \hskip 40pt v=R^{-1}Gv
$$
and $Gv$ is an eigenvector of $R^{-1}$ with same eigenvector.
It is clear that $Gv$ has to be an eigenvector of $R$ with eigenvalue
$\zeta^{-1}$. Thus $G$ maps an eigenvector corresponding to an
eigenvalue $\zeta$ to the eigenvector corresponding to its inverse.
As $R$ is diagonalizable, the geometric multiplicity of each
eigenvalue is the same as the algebraic multiplicity which is one.
The eigenvectors form a basis and thus $G$ is completely determined.
We thus have that $G^2=Id$, and up to a change of basis, $G$ is the
matrix representation of an idempotent permutation. That is, it
consists in a mere exchange of two subspaces of dimension $n/2$.
}

When the drift vector $\mathcal{S}$ is zero, the eigenvalue
$\beta$ of $M_\sigma$ is not in the spectrum of $R$.
Then the characteristic polynomial $P(x)$ of $M_\sigma$ admits
the factorization 
\begin{equation}\label{eq:ZeroDriftPolynomial}
P(x)=p(x)\tilde p(x) q(x)
\end{equation}
where the polynomials
$p$ and $\tilde p$ are irreducible of degree $n$, while
$q$ is a polynomial of degree $N-2n$.

Summing up the results of this section we have

\begin{proposition}\label{proposition:DriftConditions}
Consider a uniquely ergodic algebraic IET on $N$ intervals, 
degree $n$ and maximal rank. 
\begin{enumerate}
\item [i)] If the drift vector is non-zero and $n$ is odd,
then renormalization is not possible;
\item [ii)] if the drift vector is zero then $2n<N$.
\end{enumerate}
\end{proposition}

Note that if we have only two intervals, the two lengths are of the
form $\Omega_1$  and $1-\Omega_1$ and the condition $S=0$ would
imply $|\Omega_1| k+(1-|\Omega_1|)l=0$ for two integers $k$ and $l$
which is impossible as $|\Omega_1|$ is irrational. In the case when
$n=N=2$, the drift vector must be non-zero and yet the map
is known to be renormalizable, from the Boshernitzan and
Carroll theorem \cite{BoshernitzanCarroll}.

\subsection{Recursive tiling property}
We already introduced the standard coding, with respect to the
partition of the interval $\Omega$.
In \cite{LPVYocc,KouptsovLowensteinVivaldi,LowensteinKouptsovVivaldi}
and elsewhere, a different type of coding is also used, which is based
on the so-called {\em recursive tiling\/} property. As we mentioned
earlier (in the proof of theorem \ref{theorem:rhounit}), for each $k$,
the interval $\Omega$ splits into $\sum_{i=1}^N\sigma^k(i)$
sub-intervals, each corresponding to a $k$-cylinder in the symbolic space.
For a given point $x$, knowing for each $k$ in which cylinder $x$ lies
is the essence of this new coding, which, for reasons explained later,
is also called the Vershik coding.

We consider the substitution $\sigma$ introduced above, using
the notation $\omega^i=\sigma(i)$.
A {\it prefix\/} $\mu$ of a word $\omega_1\cdots\omega_k$ is
a sub-word of the type $\omega_1\cdots\omega_j$, with $j<k$.
Let $\mathcal{P}$ be the set of all prefixes of the words
$\omega^1,\ldots,\omega^N$, together with the $N$ additional
symbols $\epsilon_1,\ldots,\epsilon_N$, which will be regarded 
as prefixes of the 1-symbol words $\omega^1_1,\ldots,\omega^N_1$, 
respectively.
We use the notation $\mu^i\in\mathcal{P}$ to characterize
an arbitrary prefix of $\omega^i$. Associated to this notation
there is the function $\chi$ which extracts the ``exponent''
of a prefix, meaning that $\chi(\mu)=i$ if $\mu$ is a prefix
of $\omega^i$. Finally, given a prefix $\mu$ of $\omega$,
we denote by $\mu_+$ the symbol immediately following
$\mu$ in $\omega$. This quantity is well-defined, since
$\mu$ a proper sub-word of $\omega$.

If we let
$E_{\omega^i_1\ldots\omega^i_p}x=
x+\tau_{\omega^i_1}+\ldots+\tau_{\omega^i_p}$,
then for $\mu^i\in\mathcal{P}$ we have
$
E_{\mu^i}\rho \Omega_i\subset\Omega_{{\mu}^i_+},
$
and hence
$$
\bigcup_{\mu\in\mathcal{P}\atop\mu_+=i}
E_\mu\rho\Omega_{\chi(\mu)}\subset\Omega_i
$$
and
$$
\Omega=\bigcup_{i=1}^N\bigcup_{\mu\in\mathcal{P}\atop\mu_+=i}
E_\mu\rho\Omega_{\chi(\mu)}.
$$
The same decomposition can be made on each scaled sub-intervals and
we thus have a {\em recursive tiling} of $\Omega$. A given point
$x\in\Omega$ can thus be located if we know for each $k$ in which
scaled sub-interval it lies. We can show that there exists a sequence
of prefixes $(\mu_1,\mu_2,\ldots)\in\mathcal{P}$ such that
$$
x\in\bigcap_{k=1}^\infty E_{\mu_1}\rho\ldots
E_{\mu_k}\rho\Omega_{\chi(\mu_k)}.
$$
There is a consistency constraint on the sequence, namely
\begin{equation}\label{ruleg}
\chi(\mu_k)=(\mu_{k+1})_+.
\end{equation}
On the possible sequences in $\mathcal{P}^\mathbb{N}$, the dynamics
of $E$ translates into an odometer-like process which can be
conveniently represented as a Vershik map, thus justifying the
terminology of ``Vershik code''.

In \cite{LPVYocc}, the periodic sequences of the Vershik code played
an important role, firstly because it is possible to compute the
associated points explicitly, and secondly because all the
points of the lattice turned out to have periodic or 
eventually periodic Vershik codes.
By construction, a point $x$ with periodic Vershik code
$(\mu_1\ldots\mu_T)^\infty$ is an invariant point of a
contraction
\begin{equation}\label{eq:PeriodicPoint}
E_{\mu_1}\rho\ldots E_{\mu_T}\rho x=x\qquad\textrm{or}\qquad
x=\frac{1}{1-\rho^T}\sum_{i=0}^{T-1} \tau_{\mu_{i+1}} \rho^{i}
\end{equation}
where $\tau_{\mu_i}$ are the translations corresponding to the
maps $E_{\mu_i}$, which belong to $\M$. 
The rightmost expression in (\ref{eq:PeriodicPoint}) may be 
regarded as an expansion of $x$ in the base $\rho$ where
the ``digits'' $\tau_{\mu_i}$ are algebraic numbers (this 
viewpoint is developed in \cite{LowensteinVivaldi}).
Because the factors $1-\rho^T$ at denominator are not units, 
in general a periodic point is not an algebraic integer, but 
rather a point of the field $K$. 
By the decomposition (\ref{layers}), there exists
a representative $\xi\in\Xi$ such that $x$ has its orbit on
$\xi+\mathcal{M}$.
We note that the set $\xi+\mathcal{M}$ is not invariant under
scaling by $\rho$, that is,
$$
\rho\xi\in
\xi'+\mathcal{M}\qquad\textrm{where}\qquad\xi'\equiv\rho\xi\mod{\mathcal{M}}.
$$
Following \cite{LPVYocc}, we notice that, for any positive
integer $m$, the lattice $\mathcal{M}/m$ is an ${\cal O}$-module
(see remarks following equation (\ref{eq:O})),
and in particular is invariant under multiplication by $\rho$.
(Alternatively, note that multiplication by $\rho$ is represented
as an integer matrix $R$ which preserves the lattices ${\bf L}/m$.)
With our definition of $\Xi$, we then see that $\Xi\cap\mathcal{M}/m$
has $m^3$ elements, which limits the number of possible $\xi'$ after
multiplication.
Thus, there exists an integer $t$ such that
$\rho^t\xi\equiv\xi\mod{\mathcal{M}}$,
and the smallest such $t$ is called the {\it order\/} of $\xi$.
From equations (\ref{eq:PeriodicPoint}), if the point $x=\beta+\xi$ with
$\beta\in\mathcal{M}$ and $\xi\in\Xi$ has period $T$, we
find that $\xi-\rho^T\xi\in\mathcal{M}$, which
implies that $T$ must be a multiple of the order of $\xi$.

The period $T$ determines the maximum denominator of the 
fractional component $\xi$ of a periodic point, which is given
by the determinant of $I-R^T$.
Indeed, it is easy to see that the characteristic polynomial
of $I-R^T$ is $p_{T}(1-t)$ where $p_{T}$ is the
characteristic polynomial of $R^T$. Then the eigenvalues of this
matrix are of the form $1-\zeta_i^T$, where $\zeta_i$ runs through
the eigenvalue of $R$. 
From the previous section we know that we the eigenvalues of $R$
come in reciprocal pairs $(\zeta,1/\zeta)$, and a straightforward 
computation shows that
$$
d_T=|\det(I-R^T)|=\prod_{i=1}^{n/2}\frac{(\zeta^T_i-1)^2}{\zeta^T_i}
$$
where the product is taken over all the $n/2$ eigenvalues greater 
than one. There exists then a
value of the period above which $d_T$ increases monotonically with $T$.

Next we prove a bound for points with eventually periodic code.
\begin{proposition}\label{proposition:Bound}
Let $x$ be a point whose Vershik code is eventually
periodic with transient $t$ and period $T$, and let $z=\phi(x)$.
Then there exists a positive constant $C$ such that 
$$
||z||\leq C||R||^{t+T},
$$
where $||R||$ is the induced norm of matrices.
\end{proposition}
{\bf Proof:} On the corresponding lattice, we have
$$
z=\frac{1}{I-R^T}\sum_{j=0}^{T-1}R^j v_{\mu_{j+1}}
\hskip 40pt
v_{\mu_j}=\phi(\tau_{\mu_j})
$$
and we see that
$$
||z||\leq ||(I-R^T)^{-1}||\sum_{j=0}^{T-1}||R^j
v_{\mu_{j+1}}||\leq||(I-R^T)^{-1}||\max_k\{||v_{\mu_{k}}||\}\sum_{j=0}^{T-1}||R||^j,
$$
and thus
$$
||z||\leq C||(I-R^T)^{-1}||\frac{||R||^T-1}{||R||-1}.
$$
By the discussion above on the eigenvalues of $(I-R^T)$, we see that
$$
||(I-R^T)^{-1}||=\max_k\left|\frac{1}{\zeta_k^T-1}\right|=\max_k\left|\frac{1}{\zeta_k^{-T}-1}\right|
=\max_k\left|\frac{\zeta_k^T}{1-\zeta_k^T}\right|.
$$
Now the expression in the $\max$ is very small if $\zeta_k<1$ and
close to one if $\zeta_k>1$, it is clear that the maximum is
attained for the maximum value of $|\zeta_k|$ (at least for
sufficiently large $T$), which happens to be $||R||$ and
$$
||(I-R^T)^{-1}||=\max_k\left|\frac{\zeta_k^T}{1-\zeta_k^T}\right|=\frac{||R||^T}{||R||^T-1}=1+\frac{1}{||R||^T-1}.
$$
We finally have
$$
||z||\leq
\left(\frac{\max_k\{||v_{\mu_{k}}||\}}{||R||-1}\right)(||R||^T-1)\leq
C||R||^T.
$$
Let us now assume that the point $z$ has eventually periodic
Vershik code $(\mu_1\ldots\mu_t(\mu_{t+1}\ldots\mu_{t+T})^\infty)$
with transient length $t$ and period $T$.
For each positive integer $t$, we define the mapping
$$
\nabla: {\bf L} \to {\bf L}
\hskip 40pt
z\mapsto \psi_{\mu_1}R\cdots\psi_{\mu_t}R z
   =R^t z+\sum_{j=0}^{t-1}R^jv_{\mu_{j+1}}.
$$
Then $\nabla_t^{-1}x$ has $T$-periodic Vershik code, and from the
above we see that
$$
||z||\leq C||\nabla_t||||R||^T
 \leq C||R||^T(||R||^t+\sum_{j=0}^{t-1}||R||^t||v_{\mu_{j+1}}||).
$$
Applying again the geometric sum formula we see that we have a
positive constant $C'$ such that
$$
||z||\leq CC'||R||^T||R||^t.
$$
\qed\\

From this proposition, we know that the number of distinct
points $z$ with prefix length $t$ and period $T$ which belong
to a given layer $\xi+\M$ is of the order of $||R||^{(n-1)(t+T)}$. 
We also know that we have up to $d_T$ such layers compatible 
with this specific period $T$.

The number of admissible prefixes of length $t$ is given by the
admissibility matrix of the graph described by the relation
(\ref{ruleg}). This graph can be described either as a Bratteli
diagram or as a more conventional connected graph. In any way, it
can be shown \cite{GPself} that if $\sigma$ is a primitive 
substitution, then the admissibility matrix is primitive as well; 
moreover the maximum eigenvalue of this matrix is equal to the
maximum eigenvalue of $M_\sigma$, the incidence matrix of the
substitution. The number of admissible paths in the 
aforementioned graph, as well as all the number of $T$-cycles, 
are then of the order of powers of $\beta=||M_\sigma||$.

We now fix a period $T$, and consider all the prefix-lengths $t$.
We see that in order that all the points fit on all the
possible lattices, we must have
$$
d_T ||R||^{(n-1)(t+T)}\geq C ||M_\sigma||^{t+T}
$$
for some constant $C$. As $t$ increases, this is only possible 
if $||R||^{n-1}\geq ||M_\sigma||$.

We have the following result.
\begin{proposition}\label{proposition:PositiveDensity}
If
\begin{equation}\label{eq:RM}
||R||^{n-1}=||M_\sigma||,
\end{equation}
then for every period $T$ there is at least one layer
$\mathcal{M}+\xi$ with the property that the set of
(Vershik) eventually periodic points of period $T$ on
the corresponding lattice has positive density.
\end{proposition}
{\bf Proof:} The proof follows from the discussion above. 
Here we mention that the map sending each Bratteli code to a
point in $\Omega$ is not one-to-one; however, it can be shown 
\cite{GPself} that it is finitely many to one and thus
it does not change the result. \qed\\

If $\mathcal{S}\neq 0$ and $n>2$, then the identity 
(\ref{eq:RM}) is never verified. Indeed, in this case,
both matrices share the same maximum eigenvalue $\beta$,
see comments at the end of section \ref{section:algconsrenorm}.
This is consistent with the remark we made above, that an orbit
must have zero density in these cases.

In the case studied in \cite{LPVYocc}, we had
$$
||R||^2=||M_\sigma||,
$$
though the result stated here is much less powerful that the results
of \cite{LPVYocc}.

\subsection{Remarks on asymptotic behaviour}
We begin by considering the orbit of the point $z$
which corresponds to the fixed point of the substitution. We know
that there is an integer $j$ such that the code of $z$ begins with
$\sigma^k(j)$ for all $k$. In this case, we have
$$
\begin{array}{lll}
\psi^{|\sigma^k(j)|}(z)&=&z+\sum_{i=1}^N|\sigma^k(j)|_iv_i\\
&=&z+\pi\sum_{i=1}^N\langle M^k_\sigma e_j,e_i \rangle e_i\\
&=&z+\pi M_\sigma^ke_j\\
&=&z+R^kv_j.
\end{array}
$$
We clearly see that the behaviour here is determined by the largest
eigenvalue of $R$. If $\mathcal{S}\neq 0$, we know that this eigenvalue
is $\beta$, and since $|\sigma^k(j)|$ also scales like $\beta^k$,
we have a dominant linear behaviour, as expected.

If $\mathcal{S}=0$, then the behaviour is dominated by the largest
eigenvalue of $M_\sigma$ which is also part of the spectrum of $R$.
We have
$$
\pi D_{|\sigma^k(j)|}=R^k v_j=\pi M^k_\sigma e_j
$$
where the vector $D$ of discrepancies was defined in (\ref{eq:Discrepancies}).
The escape rate of the point is given by the discrepancy-vector of
the point $z$. There are known results about the discrepancy of
substitution dynamical systems, including a theorem from Adamczewski
\cite{Adamczewski} which gives the growth rate of the maximum component
of $D_k$, depending on the properties of the spectrum of $M_\sigma$.
In particular, we have an upper bound of the form $(\log k)^a k^b$, 
with $a$ a positive integer and $0<b<1$. We have this growth when the
second largest eigenvalue $\beta_2$ of $M_\sigma$ is greater than 1
(in modulus). The quantity $b$ has the form
$$
b=\frac{\log|\beta_2|}{\log\beta}
$$
and $a+1$ is the multiplicity of $\beta_2$

This indicates that a power-law behaviour can only occur when
$\beta_2$ is greater than 1 in absolute value and has 
multiplicity one. This eigenvalue
then dominates the evolution. Indeed if we have $||\pi
D_{|\sigma^k(j)|}||\sim |\beta_2|^k$, since 
$|\sigma^k(j)|\sim\beta^k$ we are led to
$$
||\pi D_{k}||\sim k^{\frac{\log|\beta_2|}{\log\beta}}.
$$
In the case studied in \cite{LPVYocc}, direct computation
gives
$$
\frac{\log|\beta_{2}|}{\log\beta}=\frac{1}{2}
$$
which leads to positive density orbits on the lattice.

In general case the computation is less straightforward, but thanks to
the recursive tiling property we can still derive estimates for the
asymptotic behaviour. Let us thus consider a point $x$ having Vershik
code $(\mu_1\ldots\mu_k\ldots)$. Then
$$
x\in E_{\mu_1}\rho\Omega_{\chi(\mu_1)}\cap E_{\mu_1}\rho
E_{\mu_2}\rho\Omega_{\chi(\mu_2)}\cap\ldots
$$
and we can deduce that (see \cite{GPself} for details)
$$
\begin{array}{lll}
x\in E_{\mu_1}\rho\Omega_{\chi(\mu_1)} & \Rightarrow &
\omega(x)=s^{|\mu_1|}\sigma(\chi(\mu_1))\\
x\in E_{\mu_1}\rho E_{\mu_2}\rho\Omega_{\chi(\mu_2)} & \Rightarrow &
\omega(x)=s^{|\mu_1|}\sigma(s^{|\mu_2|}\sigma(\chi(\mu_2))),
\end{array}
$$
where here $s$ denotes the shift map. To estimate the escape 
rate of the lattice point corresponding to $x$, we need to know
the asymptotic number of occurrences of each letter in the words above.
Let us call $L^k$ the vectors of the number of occurrences of
each of the letters $1,\ldots,N$ in the $k$-th word above
$$
L^k_i=|s^{|\mu_1|}\sigma(\cdots
s^{|\mu_k|}\sigma(\chi(\mu_k))\cdots)|_i.
$$
The number of points in the beginning of the orbit of $x$ thus
described is given by the length of the word above;
since $L^k$ is a positive vector, this is given by
$$
\sum_{i=1}^N L^k_i=||L^k||.
$$
Each shift cuts the number of every letter by a finite amount: the
number of possible shift lengths is equal to the number of possible 
prefixes, and each one decreases each component of $L^k$ in a 
finite numbers of ways. 
If we call $\eta$ the maximum length of a prefix, then
there exists not more than $N^\eta$ positive $N$-vectors
$\varepsilon_j$ of norm not exceeding $N$ such that
$$
\begin{array}{lll}
L^1&=&M_\sigma e_{\chi(\mu_1)}-\varepsilon_{\alpha_1} \\
L^2&=&M_\sigma (M_\sigma
e_{\chi(\mu_2)}-\varepsilon_{\alpha_2})-\varepsilon_{\alpha_1}\\
 &\,\vdots &\\
L^k&=&M^k_\sigma e_{\chi(\mu_k)}-\sum_{j=0}^{k-1}
M^j\varepsilon_{\alpha_{j+1}},
\end{array}
$$
where $1\leq \alpha_i\leq N^\eta$. 
As $k$ grows, the escape rate of $z$ under
$\psi$ is thus given by the asymptotic behaviour of
$$
\sum_{i=1}^N L^k_i v_i=\pi L^k=R^k v_{\chi(\mu_k)}-\sum_{j=0}^{k-1}
R^j\pi\varepsilon_{\alpha_{j+1}}.
$$
Since all integer vectors involved have uniformly bounded norm,
we have the upper bound
$$
||\psi^{||L^k||}(z)||=||\pi L^k||\leq K\sum_{j=0}^k||R||^j\leq
K\frac{||R||^{k+1}-1}{||R||-1}.
$$
This is not very sharp, yet there is no clear lower bound as one
could adjust the $\varepsilon_j$ so as to make the $L^k$ bounded.
The quantity $L^k$ is easier to express if $x$ has periodic or 
eventually periodic Vershik code. 

In what follows, we state results only for periodic codes,
to which an eventually periodic code can be reduced via
a finite number of iterations of $E$. 
In the previous section, we saw that under suitable conditions, 
these lattice points constitute a ``significantly large'' set, 
thus justifying our interest in asymptotic estimates.
We have the following result.
\begin{proposition}\label{proposition:AsymptoticPeriodic}
If a point $x$ has periodic Vershik code of period $T$, then
the corresponding lattice point $z=\phi(x)$ satisfies the 
asymptotic estimate
$$
||\psi^{kT}(z)||\sim k^{\frac{\log||R||}{\log||M_\sigma||}}
\qquad k\to\infty.
$$
\end{proposition}
{\bf Proof:} Let $x$ have periodic Vershik code $(\mu_1\ldots\mu_T)^\infty$.
By formula (\ref{eq:PeriodicPoint}), every such point is representable 
as a lattice point. 
We define
$$
\delta_T=\sum_{j=0}^{T-1} R^j\pi\varepsilon_{\alpha_{j+1}}
$$
for some indices $\alpha_j$, and 
$$
G(z)=R^T z+\delta_T.
$$
The sequence
$$
\pi L^{k T}=G^k(v_{\chi(\mu_{kT})})
$$
is an arithmetic-geometric sequence with respect to $k$. 
We use the following lemma
\begin{lemma}
If $u_k$ is an arithmetic-geometric sequence, that is
$$
u_{k+1}=a u_k+b
$$
where $a$ is a square matrix such that $I-a$ is invertible 
and $b$ is a vector, then
$$
u_k=a^k(u_0-l)-l\quad\textrm{and}\quad\sum_{j=0}^k
u_j=(I-a)^{-1}(I-a^{k+1})(u_0-l)+kl
$$
with $l=(I-a)^{-1}b$.
\end{lemma}

\noindent
{\bf Proof:} We verify that the sequence
$$
v_k=u_k-l
$$
is a geometric sequence, indeed
$$
\begin{array}{lll}
v_{k+1}&=&au_k+b-(I-a)^{-1}b\\
&=&a(u_k+a^{-1}b-a^{-1}(I-a)^{-1}b)\\
&=&a(u_k+(a^{-1}(I-a)-a^{-1})(I-a)^{-1}b)\\
&=&a(u_k-(I-a)^{-1}b).
\end{array}
$$
The result then comes by straightforward application of the classic
formulae for geometric sequences.
\qed

We now apply the lemma to our case, with $a=R^T$ and $b=\delta_T$.
Provided that $v_{\chi(\mu_{kT})}=(I-R^T)\delta_T$, we obtain
$$
||\psi^{||L^{kT}||}(z)||=||\pi L^{kT}||\sim ||R||^{kT}.
$$
Furthermore, if we define
$$
\tilde\delta_T=\sum_{j=0}^{T-1} M_\sigma^j\varepsilon_{\alpha_{j+1}}
$$
with the corresponding indices $\alpha_j$ and if we call
$$
\tilde G(z)=M^T z+\tilde\delta_T,
$$
then
$$
L^{k T}=\tilde G^k(e_{\chi(\mu_{kT})}).
$$
Applying again the argument above with $a=M_\sigma^T$ and
$b=\tilde\delta_T$ yields
$$
||L^{kT}||\sim ||M_\sigma||^{kT}
$$
which leads to
$$
||\psi^{kT}(z)||\sim k^{\frac{\log||R||}{\log||M_\sigma||}}.
$$
\qed

This estimate is very similar to the one obtained for
the orbit of zero, a point of constant Vershik code.
Moreover, as the motion on the lattice is confined to a
$(n-1)$-dimensional slab, if
$$
\frac{\log||R||}{\log\beta}=\frac{1}{n-1}
$$
we can expect positive density orbits. 

\section{Examples}\label{examples}

In section \ref{subsection:ConnectionsWithDriftVector}
we saw that an algebraic IET of even degree 
and non-zero drift vector must fulfil certain
conditions in order to be self-similar:
the contraction ratio has to be a root of an irreducible
self-reciprocal polynomial.
A question then arises, naturally, is there any such IET?
To answer this question, we shall use Rauzy-Veech induction,
a powerful tool in IET theory, described in the 
pioneering papers of G.~Rauzy and W.~Veech 
\cite{Rauzyinduce,Veech}, and in many subsequent works.

Following \cite{Veech}, we let $E$ be an $N$-interval IET on the unit interval $\Omega$, with irreducible permutation $\pi$ and lengths vector $\Lambda$.  We define
$$
I_0=[0,\sum_{i=1}^{N-1}\Lambda_i) \quad\textrm{and}\quad
I_1=[0,\sum_{i=1}^{N-1}\Lambda_{\pi^{-1}(i)}).
$$
Induction will be on the larger of these two intervals, $I$.  Computing the first-return maps, one verifies that there are two cases:

\begin{itemize}
\item{$|I_0|<|I_1|$ and the induced map has lengths vector $A_0(\pi)^{-1}\Lambda$ and permutation $a_0\pi$,}
\item{$|I_0|>|I_1|$ and the induced map has lengths vector $A_1(\pi)^{-1}\Lambda$ and permutation $a_1\pi$.}
\end{itemize}

The permutations $a_i\pi$ and the matrices $A_i$ are defined as follows:
$$
a_0\pi j=\left\lbrace
\begin{array}{ll}
\pi j & \pi j\leq\pi N\\
\pi j+1 & \pi N\leq\pi j<N\\
\pi N+1 &\pi j=N
\end{array}\right.\quad
a_1\pi j=\left\lbrace
\begin{array}{ll}
\pi j & j\leq\pi^{-1}N\\
\pi N & j=\pi^{-1}N+1\\
\pi(j-1) &\textrm{ other }j
\end{array}\right.
$$
and
$$
A_0(\pi)=\left(
\begin{array}{c|c}
\left(\begin{array}{ccc} & & \\
& {\bf 1}_{N-1} & \\
& & \end{array}\right) &
\begin{array}{c}0\\\vdots\\0\end{array}\\\hline
0\cdots 010\cdots & 1
\end{array}
\right),\quad 
A_1(\pi)=\left(
\begin{array}{c|c}
\left(\begin{array}{ccc} & & \\
& {\bf 1}_{\pi^{-1}N} & \\
& & \end{array}\right) &
\begin{array}{c}\bigcirc\\ \\
\begin{array}{ccccc}1&0&0&\cdots&0\end{array}
\end{array}\\ \hline
\bigcirc & 
\begin{array}{ccccc}
0&1&0&\cdots&0\\
0&0&1&\cdots&0\\
&&\vdots&&\\
0&0&0&\cdots&1\\
1&0&0&\cdots&0
\end{array}
\end{array}
\right)
$$
where here ${\bf 1}_n$ denotes the $n$-dimensional identity matrix and the 
first 1 in the last line of $A_0(\pi)$ occurs at the $\pi^{-1}N$ position. 
Also, if $\pi^{-1}N=N-1$, the lower-right block of $A_1(\pi)$ contains 
the single element $1$.

The construction above is naturally associated with a {\it Rauzy graph}, 
whose vertices are the irreducible permutations on $N$ symbols. 
Each vertex has two outgoing edges, labelled 0 and 1,  
ending at the vertices labelled by the permutations of the induced 
IETs on $I_0$ and $I_1$. 
Each connected component of a Rauzy graph corresponds to a 
{\it Rauzy class} of IETs. 
The Rauzy classes for 3- and 4-interval IETs are shown 
in figures \ref{classd3}--\ref{classg2}.

An IET on $N$ intervals is said to satisfy the {\it Keane property} 
\cite{Keane1} if the orbits of its discontinuity points are infinite and distinct.  
Such an IET is both minimal and uniquely ergodic;  its unique invariant 
probability measure is Lebesgue measure. 
If $E$ is an IET on $N$ intervals with irreducible permutation and 
satisfying the Keane property (this holds for almost all length vectors), 
then, according to a theorem of Rauzy \cite{Rauzyinduce}, the induced 
map described above is also an IET on $N$ intervals with irreducible 
permutation and also satisfies the Keane property.

Rauzy induction provides us with an effective tool for constructing 
self-similar IETs with a prescribed number of intervals.  
Indeed, a self-similar IET corresponds to a loop in the Rauzy graph, 
with the lengths of the induced map given by a positive eigenvector 
of the $N$-dimensional matrix obtained by taking the product of the 
$A_i(\pi)^{-1}$ around the loop, and a contraction factor given by 
the associated eigenvalue. 

\subsection{Non-zero drift vector}
\label{section:NonZeroDrift}
\begin{figure}[h]
\begin{center}
\includegraphics[width=0.7\textwidth]{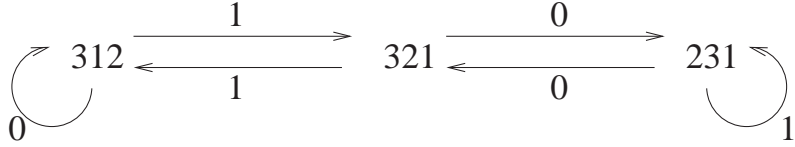}
\caption{The Rauzy diagram for 3 intervals.} \label{classd3}
\end{center}
\end{figure}

For three intervals, the Rauzy diagram has one component and
three vertices ---see figure \ref{classd3}. 
From proposition \ref{proposition:DriftConditions} $i)$, the drift
vector cannot be zero, for otherwise $n$ would have to be 1.
But then $n$ must be even, from proposition 
\ref{proposition:DriftConditions} $ii)$;
as we clearly have $n\leq 3$, we can only have quadratic cases.
This means that the characteristic polynomial $p$ cannot be irreducible.

\begin{figure}[!h]
\begin{center}
\includegraphics[width=.7\textwidth]{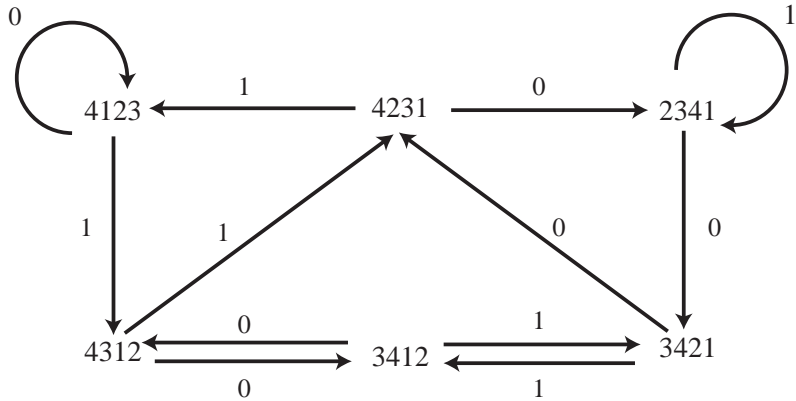}
\caption{The first Rauzy class for 4 intervals.} \label{classg1}
\end{center}
\end{figure}

\begin{figure}[!h]
\begin{center}
\includegraphics[width=\textwidth]{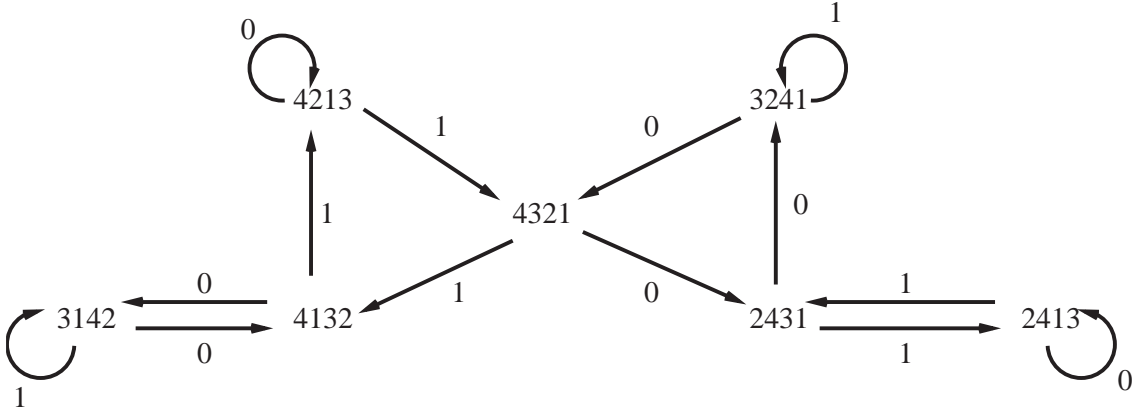}
\caption{The second Rauzy class for 4 intervals.} \label{classg2}
\end{center}
\end{figure}

For four intervals, there are two Rauzy classes, shown in figures
\ref{classg1} and \ref{classg2}.  
To find self-similar IETs,  we checked all possible cycles of length 
at most 16 in a given Rauzy class and performed the corresponding 
products of the matrices $A_i(\pi)$. 
To obtain degree-4 maps, we kept only those cycles with primitive 
product matrices and irreducible characteristic polynomials.  
For such a cycle, initiated at a specific permutation, the lengths 
of the map are the components of the normalized Perron eigenvector 
of the output matrix.
%
%
%
%

\begin{figure}[h]
\begin{center}
\includegraphics[width=\textwidth]{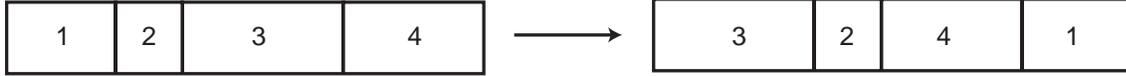}
\caption{The interval exchange transformation $E$}\label{ietpict}
\end{center}
\end{figure}

Our first example of quartic IET over four intervals comes from the
8-cycle in the graph of figure \ref{classg2}, determined by
the permutations $(4213,4213,4321,2431,3241,3241,4321,4132)$. 
As in section \ref{section:algconsrenorm}, we denote by 
$B$ the product matrix of this cycle
\beq\label{eq:Binverse}
B=\left(\begin{array}{cccc}
1& 1& 1& 1\\
0& 2& 1& 0\\
1& 2& 2& 1\\
1& 1& 1& 2\end{array}\right)^{-1}
\eeq
with characteristic polynomial
$$
p(x)=x^4-7x^3+13x^2-7x+1.
$$
We see that $p(x)$ is indeed self-reciprocal and we verify that all
his roots are real. 
The smallest eigenvalue of $B$ 
$$
\rho=\frac{1}{4}\left(7+\sqrt{5}-\sqrt{38+14\sqrt{5}}\right),
$$
corresponds to the normalized eigenvector
$$
\Lambda=\left(\rho,\,1-4 \rho+\rho^2,\,1-4\rho+5\rho^2-\rho^3,
\,-1+7\rho-6\rho^2+\rho^3\right)
$$
whose entries are positive and sum up to unity.
These data, together with the permutation $(4213)$ define
an IET $E$, with module $\M=\Z[\rho]$, for which we choose 
the canonical basis $\nu_i=\rho^{i-1}, \, i=1,\ldots,4$.
One verifies that the induced map on the first interval $\Omega_1$ 
is just the original map $E$ rescaled by a factor $\rho$.  
Figure \ref{ietpict} shows the action of the map. 

The following matrix transforms the lengths into the translations
vector $\tau=\Delta\cdot\Lambda$:
$$
\Delta=\left(\begin{array}{cccc} 0&1&1&1\\
-1&0&1&0\\
-1&-1&0&0\\
-1&0&0&0\end{array}\right).
$$
 The IET $E$ is
conjugated to its first return in the interval $\Omega_1$ via
the substitution
$$
\sigma:\left\lbrace\begin{array}{lll}
 1 & \longrightarrow & 143\\
 2 & \longrightarrow & 143223\\
 3 & \longrightarrow & 14323\\
 4 & \longrightarrow & 1443.
\end{array}\right.
$$
We note that the transpose of the inverse of $B$ in (\ref{eq:Binverse}) 
is the incidence matrix of $\sigma$. 
This identification is generally true for Rauzy induction, due to the 
identification of the $i$th substitution string with the first-return 
itinerary of the $i$th IET sub-interval. 

Finally, we verify that the drift vector is non-zero. The integer 
translation vectors $v_i=\phi(\tau_i)$ ---see equation (\ref{eq:vi})---
are
$$
(1,-1,0,0),\;(1,-5,5,-1),\;(-1,3,-1,0),\;(0,-1,0,0),
$$
so that 
$$
{\cal S}=\sum_{i=1}^4 \Lambda_i v_i=(\rho-4\rho^2-\rho^3,-1+16\rho^2-4\rho^3,4-16\rho+\rho^3,-1+4\rho-\rho^2)\
$$
is non-zero (since $\rho$ is an algebraic number of degree four).

The cycle of length 8 containing $E$ is the only one with a quartic 
irreducible polynomial. 
The other members of the cycle give 7 other IETs which renormalize 
with a scale factor $\rho$, but with different permutations and 
sets of widths in the ring $\Z[\rho]$.  
The example $E$ has the nice property that the rescaled interval is 
just $\Omega_1$, a feature not shared by all of its cycle-mates.  
For example, the cycle contains another IET with permutation (4213) 
for which the induced map again acts on the interval 
$(0,\rho), \rho\approx 0.227777$, but the latter is properly 
contained in $\Omega_1=[0,\mu), \; 
\mu=-1+7\rho-6\rho^2+\rho^3\approx 0.294963 $.  

Interestingly, for longer cycles we also have relatively few
irreducible quartic polynomials (hence rings).
The table below gives a summary of the numbers of suitable 
Rauzy cycles and associated polynomials for cycle length up to 16.  

\begin{center}
$$
\begin{array}{|c|c|c|}
\hline
\textrm{cycle length }& \sharp \textrm{ Rauzy cycles} & \sharp \textrm{ polys}\\
\hline
8 & 1 & 1\\
9 & 6 & 3\\
10 & 7 & 4\\
11 & 30 & 10 \\
12 & 80 & 27 \\
13 & 148 & 50 \\
14 & 360 & 108 \\
15 & 798 & 211 \\
16 & 1657 & 452\\
\hline
\end{array}
$$
\label{rings}
\end{center}

\subsection{Zero drift-vector}\label{subsec:zerodrift}

Here we are interested in IET belonging to the Arnoux-Rauzy family
\cite{ArnouxRauzy} to which belongs the Arnoux-Yoccoz
example \cite{Arnoux,ArnouxYoccoz} studied in great detail in
\cite{LPVYocc,LowensteinVivaldi}. Originally, this family was 
found as a way to construct a sequence of complexity $2n+1$.
The construction is the following (we keep the original 
notation from \cite{ArnouxRauzy}).

Let $\alpha,\,\beta,\,\gamma$ be three positive real numbers such that
$\alpha>\beta+\gamma$. 
We define a six-interval-exchange map $f$ over the circle of length
$\alpha+\beta+\gamma$ as follows.
First, we divide in two halves the three intervals
$I_\alpha,~I_\beta,~I_\gamma$ with respective lengths 
$\alpha,~\beta,~\gamma$ and we swap the 
halves preserving the orientation. 
Then we rotate the whole circle by $\pi$. 
It can then be shown than the induced map on 
$f(I_\alpha)$ has the same structure over three intervals
of respective lengths
$\alpha'=\alpha-\beta-\gamma,~\beta'=\beta,~\gamma'=\gamma$. We can
iterate this process: at each step we induce on the largest interval,
and at each step one of the intervals of the induced map is larger
than the sum of the other two. We can see that this process has the
same flavour as Rauzy induction. If we define the three matrices
$$
\begin{array}{ccc}
M_1=\left(\begin{array}{ccc} 1 & 1 & 1\\
0 & 1 & 0 \\
0 & 0 & 1 \end{array}\right)\quad M_2=\left(\begin{array}{ccc} 1 & 0 & 0\\
1 & 1 & 1 \\
0 & 0 & 1 \end{array}\right)\quad M_3=\left(\begin{array}{ccc} 1 & 0 & 0\\
0 & 1 & 0 \\
1 & 1 & 1 \end{array}\right)
\end{array}
$$
then it can be shown that for any $k$ the product $M_1 M^k_2 M_3$ 
occurs infinitely many times in the possible sequences of induction.
That is, for each $k$, the Perron eigenvector of each such product
will give the lengths of a renormalizable map $E_k$ over the circle.
By construction, all members of this family have zero Sah-Arnoux-Fathi 
invariant, and hence zero drift vector (proposition \ref{proposition:SAF}). 

The case $k=1$ is the Arnoux-Yoccoz case with three intervals, with
characteristic polynomial $p(x)=x^3-7x^2+5x-1$. Its unique
real root is the scaling constant $\rho=\lambda^{-3}$, where 
$\lambda$ is the so-called tribonacci number \cite[p.~233]{PytheasFogg}. 
The lengths are $(\alpha,\beta,\gamma)=(\lambda^{-2},\lambda^{-1},1)$. 
This example has very interesting properties, which include
finite decomposition.

\begin{figure}[h]
\begin{center}
\includegraphics[width=\textwidth]{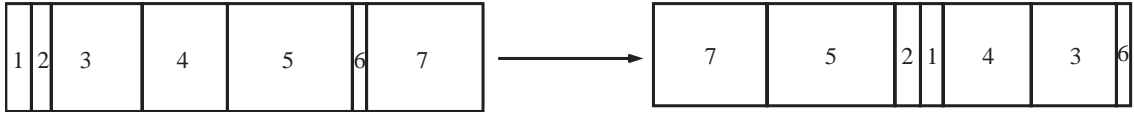}
\caption{The interval exchange $E_2$.}
\label{m2k2}
\end{center}
\end{figure}

For $k>1$, the maps maps $E_k$ are quite similar to the
Arnoux-Yoccoz map. We will present some computations 
for the cases $k=2,3,4,5,6$ (see figure \ref{m2k2} for $k=2$),
which show that they are likely to share the finite decomposition 
property. 

The module of $E_k$ is the half-integer lattices
$\M=\mathbb{Z}[\lambda_k]/2$, where $\lambda_k$ is the 
smallest real roots of the polynomial $f_k$, given by
$$
f_k(x)=x^3-(k+4) x^2+(3k+4)x -1, \quad k\geq 1.
$$
One verifies that for all $k\geq 1$, the polynomial $f_k$ has no integer 
roots and hence is irreducible.
The (non-normalized) vector of lengths of the partition is
$$
\Lambda=(\lambda_k- \lambda_k^2/2, \lambda_k- \lambda_k^2/2,1/2-3\lambda_k/2+\lambda_k^2/2,1/2-3\lambda_k/2+\lambda_k^2/2,1/2,\lambda_k/2,1/2-\lambda_k/2),
$$
and the corresponding translation vectors are
$$
\frac{1}{2}(2 + \lambda_k - \lambda_k^2,
    2 - 3\lambda_k + \lambda_k^2,
    3 - 4\lambda_k + \lambda_k^2,
    1 + 2\lambda_k - \lambda_k^2,
    \lambda_k-1, 1 - \lambda_k,
    \lambda_k-3).
$$
One can then verify directly that the drift vector is indeed zero. 
The maps $E_k$ are not renormalizable but in each cases ---as in
\cite{LPVYocc}--- the induced map in the first interval
leads to a renormalizable IET on seven intervals.
We shall carry out the computations directly with the maps $E_k$.

In order to investigate the finite decomposition property, we
will perform a fairly straightforward algorithm implemented in C.
For each $k$ we will first select the slab of the 3-dimensional
half-integer lattice consisting of the points belonging to the
interval $[0,2 - \lambda_k)$ which lie inside a cube of prescribed
size. Then we iterate forward and backward all the points in
a small cube around the origin to check how many points
in the big cube slab are attained.
This method is fast, easy to implement, and applicable to
any algebraic IET. It is however memory-intensive, which
constrained the maximum linear size of the cube we were 
able to consider to about 600 points.
For the sake of speed, when mapping a lattice point into 
the intervals we use standard double precision 
(as opposed to arbitrary precision).
It should be noted however that there is no error propagation 
in the dynamics, as each point is represented in vector 
form using the $\Z$-module basis.

To represent lattice points, we chose the 
basis $(1,\lambda_k,\lambda_k^2)$  for $\Z[\lambda]$.
It should be noted that this choice is not necessarily
the best in all cases. Indeed, orbits are often highly 
``anisotropic'', escaping much faster in a specific 
direction than in others.
As a consequence, covering a significant part of the cube 
around the origin could be time-consuming. 
While more favourable bases have been found in special cases, 
no general optimal family has been found.
For instance, for $k=1$ we have the same case as in \cite{LPVYocc}
and yet the orbits are not distributed as uniformly as in this
paper. In fact, $\lambda_1$ is the third-root of the tribonacci
number used as a basis in \cite{LPVYocc}, switching to this basis
gives more uniform orbits. This process does not generalize to 
greater $k$ though.

Nevertheless, on all the cases we have investigated we see that a
very small set of initial conditions is required to cover at least a
significant part of the lattice around zero. We define
$$
\mathcal{C}_D=\Big\{\frac{1}{2}(p,q,r)\in\frac{1}{2}\mathbb{Z}^3\,\,:\,\,
  \frac{1}{2}(p+q\lambda_k+r\lambda^2_k)\in[0,2-\lambda_k),
\,\,\, \max\{|p|,|q|,|r|\}\leq
2D\Big\}.
$$
In table \ref{calculs}, we summarize the results of the
investigations. For each $k$, we give the maximum size $D$ of
$\mathcal{C}_D$ that we were able to fill with backward and forward
iterates of points in an initial region $\mathcal{C}_d$  of size
$d$. We also give an upper bound on the number $T$ of iterations
required.

\begin{table}[h]
$$
\begin{array}{|c|c|c|c|}
\hline
k & D & d & T\\
\hline
2 & 400 & 1 & 3\times 10^6\\
3 & 400 & 3 & 6\times 10^7\\
4 & 400 & 1 & 2.6\times 10^8 \\
5 & 400 & 1 & 9\times 10^8\\
6 & 200 & 3 & 1\times 10^9\\
 \hline
\end{array}
$$
\caption{Data for lattice-filling experiments.}
\label{calculs}
\end{table}

For each $k$, the substitutions $\sigma_k$ are explicitly computable
and we can use propositions \ref{proposition:PositiveDensity} and
\ref{proposition:AsymptoticPeriodic} to give us additional evidence regarding 
the finite-decomposition hypothesis. Indeed the latter
proposition gives us asymptotic estimates on the escape rate of
orbits of the lattices having periodic Vershik code, while the
former gives a condition for these orbits to be fairly common. 
Both proposition rely on the value of the following quantity
$$
\upsilon=\frac{\log||R||}{\log||M_{\sigma_k}||}
$$
which gives the exponent of the escape rate of the orbit of a point
having periodic Vershik code. These points have density one in at
least one coset $\xi+\mathcal{M}$ if $\upsilon=1/(n-1)$. In the
following table we give approximations of the values of $\upsilon$
for several values of $k$
$$
\begin{array}{|c|c|c|c|c|c|c|c|}
\hline
k & 1 & 2 & 3 & 4 & 5 & 6 & 7\\
\hline \upsilon& 0.5 & 0.5 & 0.5 & 0.546715 & 0.595958 & 0.623202 & 0.642502\\
\hline
\end{array}
$$

We see that for $k=1,2,3$, we can deduce from the proposition
mentioned above that the points having periodic Vershik code (a set
which have positive density) have positive density, which is another
indication that $E_2,E_3$ might have the finite decomposition
property. The data for the other values of $k$ are inconclusive
though. The following table shows values of $\upsilon$ for increasing
values of $k$. The value of $\upsilon$ appears to increases slowly 
toward unity.
$$
\begin{array}{|c|c|c|c|c|c|c|c|}
\hline
k & 17 & 30 & 50 & 100 & 150 & 200 & 500\\
\hline
\upsilon& 0.721655 & 0.756235 & 0.780871 & 0.8074 & 0.820165 & 0.828247 & 0.849766\\
\hline
\end{array}
$$

\begin{figure}[h]
\begin{center}
\includegraphics[width=0.7\textwidth]{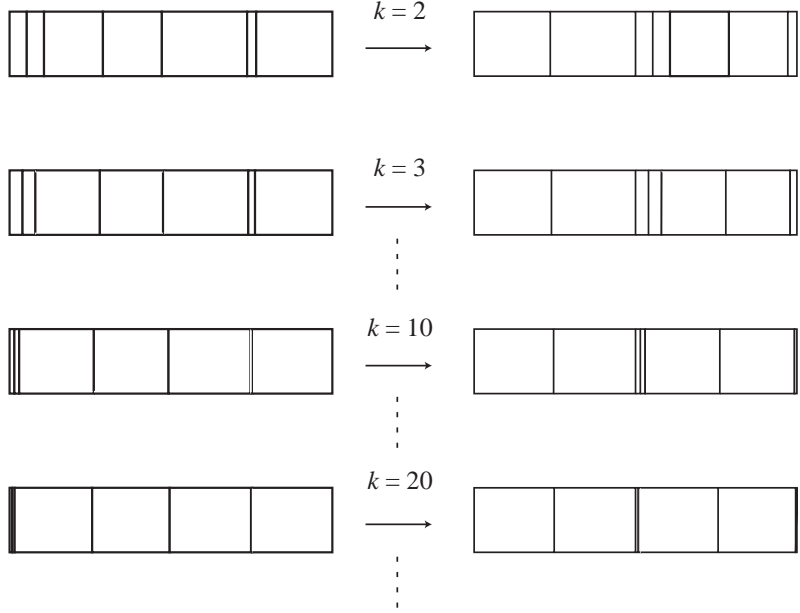}
\caption{The interval exchanges $E_k$ for increasing $k$, showing
that three intervals tends to vanish, and the dynamics approaches
a rotation with flips.} \label{increasingk}
\end{center}
\end{figure}

This can be explained if we look at the partitions of the maps $E_k$
for increasing $k$. 
As displayed in figure (\ref{increasingk}), we see that the 
intervals carrying the ``rotation-flip'' dynamics tends to fill 
the whole interval, leading to a global dynamics approaching  
a simple non-zero drift dynamics, for which $\upsilon=1$.

\subsection{Self-similar $E_2$ model}

We conclude our discussion of zero-drift-vector examples with a more 
detailed look at a self-similar version of the map $E_2$ introduced 
in section \ref{subsec:zerodrift}.  
The particular IET, which we shall call $E_2^*$, is a member of a 
29-vertex cycle in the same 294-vertex Rauzy class of 7-interval 
IETs as the self-similar Arnouz-Yoccoz model of \cite{LPVYocc}.  
Its permutation  is $(5462731)$, and the 7-dimensional cyclic 
product matrix is
$$
B=\left(\begin{array}{ccccccc}
4&9&6&6&4&8&2\\
0&2&1&1&0&1&0\\
0&2&3&0&2&2&0\\
1&2&1&2&0&1&0\\
1&1&1&1&2&2&1\\
0&2&2&0&2&3&0\\
1&1&1&1&1&1&1
\end{array}\right)^{-1}
$$
with characteristic polynomial (cf.~equation (\ref{eq:ZeroDriftPolynomial})))
$$
P(x)=(x^3-6 x^2+10 x-1)
(x^3-10 x^2+ 6 x-1)
(1-x).
$$
The scaling constant $\rho\approx 0.106711$ is the real root of the polynomial
$$
p(x)=x^3-6x^2+10x-1.
$$
One verifies that $\beta=\rho^{-1}$, the real root of 
$\tilde p(x)=x^3-10 x^2+ 6 x-1$, is a Pisot number.

The 29 members of the Rauzy cycle containing $E_2^*$ are all 
self-similar IETs with the same scale factor $\rho$.  
The map $E_2^*$ was selected for special consideration because of 
several additional simplifying features. 
In particular, its translation module is $\M=\Z[\rho]$, its 
rightmost interval $\Omega_7$ is of length $\rho$, and the 
induced IET on $\Omega_7$ is just the original map rescaled by 
$\rho$. 

\begin{figure}[h]
\begin{center}
\includegraphics[width=\textwidth]{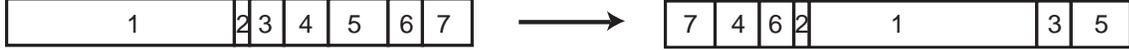}
\caption{The self-similar $E_2$ map.} \label{fig:SSE2map}
\end{center}
\end{figure}

The vector $\Lambda$ of interval lengths is the positive 
7-vector associated with $\rho$ in the diagonalization 
of $B^{-1}$, normalized to total length unity. Specifically,
choosing the canonical basis $(1,\rho,\rho^2)$ for
$\M=\Z[\rho]$, we find
$$
\Lambda=(1 - 5 \rho + 2 \rho^2, -1 + 10 \rho - 3 \rho^2, 
   1 - 9 \rho + 3 \rho^2, \rho -
   \rho^2, -1 + 11 \rho - 4 \rho^2, 1 - 9 \rho + 3 \rho^2, \rho)
$$
with corresponding integers translations $v_i=\phi(\tau_i),\,\,i=1,\ldots,7$,
$$
(v_i)=\left((0,3,-1),(0,-2,0),(1,-7,2),(-1,5,-2),(1,-8,3),(0,-6,2),(-1,1,0)\right).
$$
The IET is illustrated in figure \ref{fig:SSE2map}.
Calculating $\sum_i \Lambda_i v_i$, we verify directly that the drift 
vector vanishes. 

Iterating the map $E_2^*$, we can examine the first-return orbits of the 
7 intervals of the induced map on $\Omega_7$ and establish the 
recursive tiling property, with associated substitution rule 
$$
\sigma:\left\lbrace\begin{array}{lll}
 1 & \longrightarrow & 7114115\\
 2 & \longrightarrow & 7114121361361214115\\
 3 & \longrightarrow & 711412136136135\\
 4 & \longrightarrow & 71141214115\\
 5 & \longrightarrow & 71156136135\\
 6 & \longrightarrow & 711561361361214115\\
 7 & \longrightarrow & 7115.
\end{array}\right.
$$

We obtain the following result:

\begin{proposition}\label{proposition:E2}
We consider the mapping $E_2^*$. The following holds
\begin{itemize}
\item [$i)$]
Let $x$ be a point of the unit interval and let $\mu(x)$ be its Vershik code. 
Then $x\in\Q(\rho)\cap[0,1)$ if and only if $\mu(x)$ is eventually periodic.
\item [$ii)$] The restriction of the module $\M$ to the unit interval
consists of the (forward and backward) orbits of the discontinuity points.
\item [$iii)$]
The map has the finite decomposition property. Specifically, 
there is a positive integer $M$ such that for each $\xi\in K$ 
the restriction of the layer $\xi +\M$ to the unit interval 
decomposes into at most $M$ orbits.
\end{itemize}
\end{proposition}

We shall not provide a detailed proof of the proposition above, since 
it is a straightforward adaptation to the present model of the proof 
provided in section 3 of \cite{LPVYocc} for the Arnoux-Yoccoz cubic model. 
Let us merely summarize the main ideas.  
We represent a point  $x$ in $\Q(\rho)\cap [0,1)$, as $x=\xi+z$,
with $\xi\in\Xi$ ---see equation (\ref{eq:Xi})---
and $z\in\Z[\rho]$. Specifically, using the
basis $(1,\rho,\rho^2)$ for $\Q(\rho)$ we have
\begin{eqnarray*}
\xi&=& r_0+r_1\rho^2+r_2\rho^2\qquad r_i\in\Q\cap[0,1)\\
z &=&  m_0+m_1\rho+m_2\rho^2\qquad m_i\in\Z.
\end{eqnarray*}
The restriction of $x$ to the unit interval leads to the condition
$$
m_0=-\lfloor \xi + m_1\rho+m_2\rho^2\rfloor.
$$
Any point on the forward or backward orbit of $x$ can be represented in 
the same way, with the same $\xi$. Since each point is completely 
determined by an integer pair $(m_1,m_2)$, it is natural to associate 
the orbit with a walk on the lattice  ${\bf L}_\xi'\cong\Z^2$
(with reference to equation (\ref{eq:Isomorphism}), here we have $n=3$ and $b=1$). 

To prove the finite decomposition property of proposition 
\ref{proposition:E2}, we must show that the number of 
lattice orbits in ${\bf L}_\xi'$ is bounded uniformly in $\xi$.  
The strategy is to show that every lattice orbit visits a finite
$\xi$-independent core region containing no more than a 
finite number $M$ of lattice points. 
The existence of such a region is also the key to proving the
second statement in proposition \ref{proposition:E2}: 
one needs only to verify that the lattice orbits of the discontinuity 
points visit all of the points of the core region of ${\bf L}_0'$.

To complete the proof, one shows that for each 
$\xi\in \Xi$, some $k$th iterate of the left-shift map 
for the Vershik code is represented in ${\bf L}_\xi'$ by a map $\gamma^k$ with the 
following property: every $(m_1,m_2)\in {\bf L}_\xi'$ has a forward 
$\gamma^k$ orbit which is eventually trapped in a finite, 
$\xi$-independent core region.  This orbit must eventually settle down 
in a limit cycle of some finite period $t$. Thus the original point 
$x$ has an eventually periodic Vershik code of period $kt$, 
thereby establishing the ``only if" part of proposition \ref{proposition:E2} $i)$; 
while the ``if" part follows at once from (\ref{eq:PeriodicPoint}).
Since the periodic tail of this code coincides with that of a point 
in the limit cycle in the core region, the two points lie on the 
same $E_2^*$ orbit, and we have the remaining element in the proof 
of the finite decomposition property as well. 

\begin{figure}[h]
\begin{center}
\includegraphics[width=0.9\textwidth]{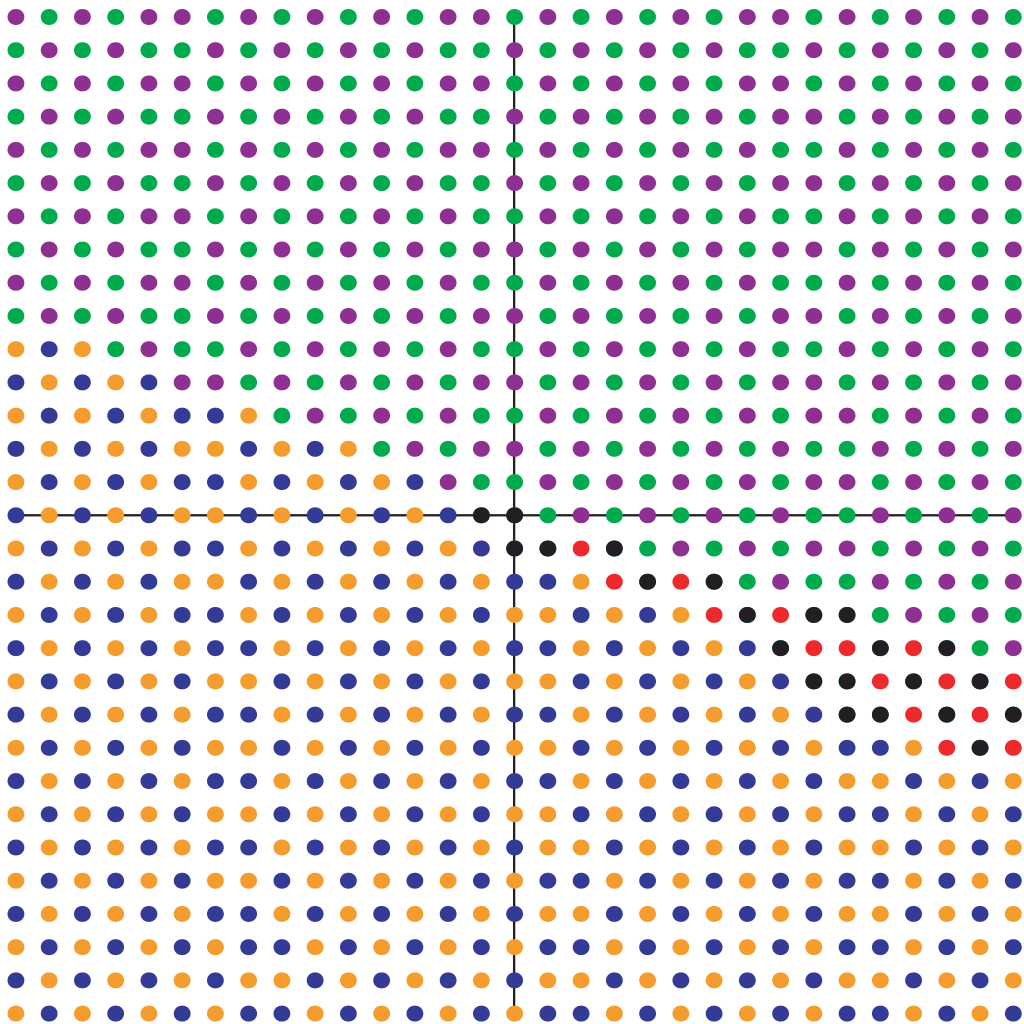}
\caption{Points of the six lattice orbits of the discontinuity points 
of the IET $E_2^*$. 
Each of the orbits is arbitrarily assigned a colour.} \label{fig:SSE2lattice}
\end{center}
\end{figure}

It is important to point out that the proof above relies 
essentially on the self-similarity of the IET, and on the
Pisot property of the reciprocal of the scale factor.
In particular, it is the latter which ensures a contraction 
to the core region.
In \cite{LPVYocc}, additional properties of the specific IET were used 
to characterize in some detail the lattice orbits in the various ${\bf L}_\xi'$. 
The orbits of the self-similar Arnoux-Yoccoz model were found to 
be non-crossing, filling out with unit density disjoint sectors 
of the lattices.  
It turns out that the lattice orbits of the $E_2^*$ model do not have 
these features. 
This is already true for the discontinuity-point orbits, for which the 
lattice is partitioned into three sectors, each filled out by a pair 
of mutually crossing orbits, neither of which is locally lattice-filling.  
The situation near the origin is shown in figure \ref{fig:SSE2lattice}.


\end{document}